\documentclass[12pt]{amsart}
\usepackage{amssymb}
\usepackage{amscd}
\usepackage{enumerate}
\usepackage{amsfonts}
\usepackage{pb-diagram,lamsarrow,pb-lams}
\usepackage[sorted]{amsrefs}
\newtheorem{thm}{Theorem}[section] \newtheorem{cor}[thm]{Corollary}
 \newtheorem{prop}[thm]{Proposition}
\newtheorem{claim}[thm]{Claim} \newtheorem{fact}[thm]{Fact}
\theoremstyle{definition}
\newtheorem{defn}[thm]{Definition}
\newtheorem{ntn}[thm]{Notation}
\theoremstyle{remark}

\newtheorem{rem}[thm]{Remark}
 \numberwithin{equation}{section}
    {\medskip\begingroup\leftskip 0.5cm\rightskip 0.5cm\noindent\begin{small}{\bf Remark.}}
    {\end{small}\par\endgroup}
\newcommand{\mrmk}[1]
{{\tiny$^{\spadesuit}$}\marginpar{\fbox{\footnotesize #1}}}
\def\strutdepth{\dp\strutbox}%
\def\marginalnote#1{\strut\vadjust{\kern-\strutdepth\specialnote{#1}}}%
\def\specialnote#1{\vtop to \strutdepth{\baselineskip%
\strutdepth\vss\llap{\hbox{\scriptsize \bf #1}}\null}}%

\newlength{\ngnut}
\newlength{\ngnutt}
\newlength{\ngnuttt}
\newcommand{\notation}[3]{
\settowidth{\ngnut}{#1}
\addtolength{\ngnut}{10pt}
\settowidth{\ngnuttt}{ Page \pageref{#3}}
\setlength{\ngnutt}{\textwidth}
\addtolength{\ngnutt}{-\ngnut}
\addtolength{\ngnutt}{-\ngnuttt}
\addtolength{\ngnuttt}{-5pt}
\begin{minipage}[t]{\ngnut} #1
\end{minipage}
\begin{minipage}[t]{\ngnutt}
#2
\end{minipage}
\begin{minipage}[c]{\ngnuttt} Page \pageref{#3}
\end{minipage}
}

\newcommand{\mfbox}[1]
{%
\medskip
\begin{center}\setlength{\fboxrule}{1pt}\setlength{\fboxsep}{10pt}
\fbox{ \begin{minipage}{0.9\textwidth} #1
\end{minipage}
}
\end{center}%
\medskip
}

\newcommand{\vlabel}{\label}
\newcommand{\elabel}[1]{\marginalnote{\boxed{#1}}\label{#1} }    

\def\fF{\mathfrak{F}}  

\def\CT{\mathcal{T}}   
\def\Mi{\mathfrak{M}}  
\def\Mib{\overline{\Mi}}

\def\Sym{Sym}
\newcommand{\C}{\mathbb{C}}
\newcommand{\Z}{\mathbb{Z}}
\newcommand{\CR}{\mathcal R} 
\newcommand{\N}{\mathbb N}

\newcommand{\rest}{\upharpoonright}
\newcommand{\co}{\circ}
\newcommand{\id}{\mathrm{id}}
\newcommand{\la}{\langle}
\newcommand{\ra}{\rangle}
\newcommand{\PP}{\mathbb{P}} 
\newcommand{\R}{\mathbb{R}}
\newcommand{\Ran}{\R_{an}}
\newcommand{\Rae}{\R_{an,exp}}

\newcommand{\sub}{\subseteq}

\def\FF{\mathbb F}

\def\CH{\mathcal{H}}

\def\CE{\mathcal{E}}
\reversemarginpar

\def\wt#1{\widetilde{#1}}

\newcommand{\leftexp}[2]{{\vphantom{#2}}^{#1}{\hspace{-1pt}#2}}
\def\tra#1{\leftexp{t}{#1}}
\def\Q{\mathbb{Q}}

\newcommand{\teta}[2]{\vartheta\!\!\begin{bmatrix}
                                   #1  \\
                                   #2
                                 \end{bmatrix}\!\!
                               }

\newcommand{\eexp}[1]{\exp\bigl( #1 \bigr)}
\newcommand{\snorm}[1]{\| #1 \|_s}
\newcommand{\pav}[2]{\CE^{#1}_{#2}}
\def\GDD{G_D(D)_0}
\def\CX{\mathcal X}

\def\spc#1{{\mathbb R}_{#1}}   
\def\e{\mathbf{e}}

\title[Theta functions]{Definability of restricted Theta functions and families of abelian varieties}
\author{Ya'acov Peterzil}
\address{University of Haifa}
\email{kobi@math.haifa.ac.il}
\author{Sergei Starchenko}
\address{University of Notre Dame}
\email{sstarche@nd.edu}
\thanks{}

\begin{document}

\maketitle

\begin{abstract} We consider some classical maps from the  theory of
abelian varieties and their moduli spaces and prove their definability, on
restricted domains, in the o-minimal structure $\Rae$. In particular, we prove that
the embedding of moduli space of principally polarized ableian varierty,
$Sp(2g,\Z)\backslash \CH_g$, is definable in $\Rae$, when restricted to Siegel's
fundamental set $\fF_g$. We also prove the definability, on appropriate domains, of
embeddings of families of abelian varieties into projective space.
\end{abstract}

\section{Introduction}

The  goal of this paper is to establish a link between the classical analytic theory
of complex abelian varieties and theory of o-minimal structures. In \cite{wp}, we
established a similar link in the one-dimensional case, by showing that the analytic
$j$-invariant and more generally, the Weierstrass $\wp$-functions are definable,
when restricted to appropriate domain, in the o-minimal structure $\Rae$ (this is
the expansion of the real field by the real exponential function as well as all
restrictions of real analytic functions to compact rectangular boxes $B\sub \R^n$.

 In  \cite{Pila}, Pila makes use of this o-minimal link in order to prove
some open cases of the Andr\'e-Oort Conjecture by applying his powerful theorem with
Wilkie \cite{Pila-Wilkie} on the connections between arithmetic and definable sets
in o-minimal structures. The results in our present paper suggest the possibility of
similar applications of tools coming from the theory of o-minimality in order to
attack more open cases of the Andr\'e-Oort  and other related conjectures.

O-minimality offers  some other powerful machinery,  and at the end of the paper we
give an example of one such application of our results to the theory of
compactifications in complex analytic spaces.

\medskip

We review briefly the basic setting in which we work here. We
consider, for every $g\geq 1$, the family of complex
$g$-dimensional tori and in particular the sub-family of abelian
varieties, namely those tori which are isomorphic to projective
varieties. Every abelian variety of dimension $g$ is isomorphic to
$\CE^D_{\tau}=\C^g/(\tau\Z^g+D\Z^g)$, where $\tau$ is a symmetric
$g\times g$ complex matrix with a positive definite imaginary part
and $D=diag(d_1,\ldots, d_g)$ with $d_1|d_2\cdots |d_g$
non-negative integers (called a polarization type of
$\CE^D_{\tau}$). Thus, for every such matrix $D$, the family of
all corresponding abelian varieties is parameterized by the {\em
Siegel half space}
$$\CH_g=\{\tau \in Sym(g,\C):
Im(\tau)>0\}.$$  For each $\tau \in \CH_g$, there is a holomorphic map
$h_{\tau}:\C^g\to \PP^k(\C)$ (for some $k$ depending only on $D$), which is
invariant under the lattice $\tau\Z^g+D\Z^g$ and induces an embedding of
$\CE^D_{\tau}$ into $\PP^k(\C)$. While $h_{\tau}$ is periodic and therefore cannot
be definable in an o-minimal structure, its restriction to a fundamental
parallelogram $E^D_{\tau}\sub \C^g$ is immediately seen to be definable in the
o-minimal structure $\R_{an}$ (since the real and imaginary parts of this map are
real analytic).

 Each $\CE^D_{\tau}$ is naturally endowed with a
polarization of type $D$.
 For each such $D$ there is a corresponding discrete group
$G_D\sub Sp(2g,\Q)$ acting on $\CH_g$ and two polarized abelian
varieties in the family are isomorphic if and only if their
corresponding $\tau$'s are in the same $G_D$-orbit.

 Let us
consider the family of principally polarized varieties, namely the
case $D=I$, in which case $G_D=Sp(2g,\Z)$.
 By \cite{baily-borel}, there is a holomorphic map $F:\CH_g\to
\PP^m(\C)$, for some $m$, inducing an embedding of $Sp(2g,\Z)\backslash \CH_g$ into
$\PP^m(\C)$, whose image is Zariski open in some projective variety. Because of the
periodicity of this map, one cannot hope that the map $F$ will be definable in an
o-minimal structure. Instead, we consider a restriction of $F$ to a subset of
$\CH_g$. It follows from Siegel's reduction theory that there is a semialgebraic set
$\fF_g\sub \CH_g$, called \emph{the Siegel fundamental set}, which contains a single
representative for every $Sp(2g,\Z)$-orbit.
 We prove here (see Theorem \ref{ppp}):
\begin{thm}\vlabel{mainn1} There is an open set $U\sub \CH_g$ containing the Siegel fundamental set,
such that the  restriction of the above $F:\CH_g\to \PP^m(\C)$ to $U$ is definable
in the o-minimal structure $\R_{an,exp}$.
\end{thm}

Note that since $F$ is $Sp(2g,\Z)$-invariant, and $\fF_g$ is a
fundamental set for $Sp(2g,\Z)$, we have $F(U)=F(\CH_g)$ and in
this sense the above theorem says that the quotient
$Sp(2g,\Z)\backslash \CH_g$ can be definably embedded into
projective space, in the structure $\Rae$.

While the map $\tau \mapsto F(\tau)$ is  a function on $\CH_g$, the theorem above
follows from the definability of a map in two sets of variables $(z,\tau)$. Namely,
for every fixed polarization type $D$, we consider the family of abelian varieties
of polarization type $D$, $\mathcal F^D=\{\CE^D_{\tau}:\tau\in \CH_g\}$. As was
noted in \cite{wp}, the family $\mathcal F^D$ can be viewed as a semi-algebraic
family of complex tori, by identifying the underlying set of each $\CE^D_{\tau}$
with its fundamental parallelogram $E^D_{\tau}$. We let $h^D_{\tau}:E^D_{\tau}\to
\PP^k(\C)$ be the embedding of the complex torus $\CE^D_{\tau}$ into projective
space and consider the family $\{h^D_{\tau}:\tau\in \CH_g\}$. Again, for reasons of
periodicity we cannot hope for this whole family to be definable in an o-minimal
structure. However, using Siegel's reduction theory we can find a fundamental set
$\fF_g^D\sub \CH_g$, such that every $\CE_{\tau}^D$ is isomorphic, as a polarized
abelian variety of type $D$, to one of the form $\CE^D_{\tau'}$, for $\tau'\in
\fF_g^D$. We prove (see Theorem \ref{Gdef2}):

\begin{thm}\vlabel{mainn2} For every polarization type $D$, there is in $\Rae$ a definable open set $U\sub \CH_g$ containing
$\fF^D_g$ and a definable family of maps $\{h^D_{\tau}:\tau\in U\}$, such that each
$h^D_{\tau}:E^D_{\tau}\to \PP^k(\C)$ is an embedding of the abelian variety
$\CE^D_{\tau}$ into projective space.
\end{thm}

 The maps $F$ and the family $\{h^D_{\tau}\}$ are given in coordinates by modular
 forms, which themselves can be described using Riemann theta functions, of the form $\teta{a}{b}(z,\tau)$, for $a,b\in \R^g$.
 Most of our work in the paper goes towards proving:
 \begin{thm}\vlabel{main3} Fix a polarization type $D=diag(d_1,\ldots,d_g)$. For each $a,b\in \R^g$, the restriction of $\teta{a}{b}$
 to the set $$\{(z,\tau)\in \C^g\times \CH_g: \tau\in \fF_g\,\, , \, z\in E^D_{\tau}\}$$ is definable in
 $\R_{an,exp}$.
 \end{thm}

  Much attention has been given classically to the possible compactifications of the quotient
  $G_D\backslash {\CH_g}$. At the end of the paper we show how removal of singularities
  in o-minimal structures can be used to show that the closure in projective space of the image
  of some of the above maps,  is an
  algebraic variety.
\\

For background on o-minimal structures, see \cite{vddm}. We should
note that the o-minimality and main properties of $\R_{an,exp}$
were established in
 \cite{vdd-miller} and \cite{vdd-M-M}. Actually, we need almost no background from
 the theory of o-minimal structures, except for the simple fact that for every
 holomorphic function $f:U\to \C$, on some open $U\sub \C^n$, if $S\sub U$ is a
 compact semialgerbaic set then $f\rest S$ is definable in the o-minimal structure $\Ran$ (so in
 particular in $\Rae$).

\section{Conventions and Notations}
 Because of the need for  heavy notation we begin with a
 list of symbols to be used in the paper.\subsection{Conventions}
For a ring $R$ we will denote by $M_g(R)$ the set of  of all
$g\times g$ matrices over $R$,
and by $Sym(g,R)$ the set of all symmetric $g\times g$ matrices over $R$.

 We always consider $a\in R^g$ as a column vector.
For $a\in R^g$ and $M\in M_g(R)$, we  denote by $\tra a$ and $\tra M$ the transpose of $a$ and $M$, respectively. \\
For a matrix $y\in M_g(R)$ and $a\in R^g$ we  denote $\tra a y
a$ by $y[a]$.

We denote complex exponentiation by $exp$.

\subsection{Notations used in the paper}\ \\
\noindent \notation{$\CH_g$}{Siegel upper half-space of degree $g$}{note:siegel}
\notation{$\Lambda^D_\tau$}{ The lattice $D\Z^g+\tau\Z^g$}{note:polar}
\notation{$\pav D \tau$}{ The torus $\C^g/\Lambda^D_\tau$}{note:polar}
\notation{$Sp(2g,K)$}{The symplectic group}{note:sp} \notation{$G_D$}{A special
subgroup of $Sp(2g,\Q)$}{note:Gd} \notation{$\Mi_g$}{The set of Minkowski reduced
matrices}{def:minred} \notation{$\fF_g(G)$}{A fundamental set for the action of
$G$}{note:fFD} \notation{$\vartheta(z,\tau), \teta a b (z,\tau)$} {Riemann Theta
functions}{note:theta} \notation{$\varphi^D$}{A holomorphic map from $\C^g\times
\CH_g$ into $\PP_n(\C)$ \\ whose components are theta functions}{note:vd}
\notation{$\varphi^D_\tau$}{The map $z\mapsto \varphi^D(z,\tau)$} {note:vdt}
\notation{$\spc n$} {The standard polyhedral cone $(\R^{\geq 0})^n$}{note:spc}
\notation{$\CT(S)$}
{A tube in $\C^n$ corresponding to a set $S\subseteq \R^n$}{defpseta}\\
\notation{$\e(z)$} {The complex function $\e(z)=e^{2\pi i z}$}{note:be} \notation{
$\Mib_g$}{An integral polyhedral cone containing $\Mi_g$}{note:miclosed}
\notation{$C_m$}{An integral polyhedral cone in $\R^g\times Sym(g,\R)$}{note:notec}
\notation{$E_\tau^D$}{The fundamental parallelogram of $\Lambda_\tau^D$}
{note:edtou} \notation{$\bar E_\tau^D$}{The topological closure of $\bar E_\tau^D$}
{note:edtou} \notation{$\CX^D(V)$ and $\CX^D_{<K}(V)$}{Subsets of $\C^g\times
\CH_g$}{note:cxdv} \notation{$\GDD$}{A special subgroup of $G_D$}{GDD}
\notation{$\fF_g^D$}{The fundamental set $\fF_g(\GDD)$}{note:fFgd}
\notation{$\Psi^D(\tau)$}{ The map $\Psi^D(\tau)= \varphi^D(0,\tau)$}{note:phipsi}
\notation{$\Phi^D(z,\tau)$}{ The map
$\Phi^D(z,\tau)=(\varphi_D(z,\tau),\Psi^D(\tau))$} {note:phipsi}

\section{Polarized abelian varieties and the Siegel fundamental domain}
In this section we review briefly  some known facts about tori and polarized abelian
varieties. We refer to  \cite{BL} for more details.

\subsection{Complex $g$-tori, abelian varieties and polarization}
For a positive $g\in \N$, by {\em a complex $g$-torus} we mean the
quotient group $\C^g/\Lambda$, where $\Lambda\subseteq \C^g$ is a
lattice, i.e. a subgroup of $(\C^g,+)$ generated by $2g$ vectors
which are $\R$-linearly independent. Under the induced structure,
a complex $g$-torus $\C^g/\Lambda$ is a compact complex Lie group
of dimension $g$, and vice versa, every $g$-dimensional compact
complex Lie group is bi-holomorphic to a complex $g$-torus.

If $f\colon \C^g/\Lambda\to \C^g/\Lambda'$ is a bi-holomorphism between two tori then
$f$ can be lifted to a $\C$-linear map $F\colon \C^g\to \C^g$ with $F(\Lambda)=\Lambda'$,
and vice-versa, every such $F$ induces an isomorphism between corresponding tori.

A torus $\C^g/\Lambda$ is called \emph{an abelian variety} if it is biholomorphic
with a projective variety in $\PP^k(\C)$ for some $k$. If $g>1$ then not every
$g$-torus is an abelian variety. The following criterion is the well-known Riemann
condition.
\begin{thm}[Riemann condition]  A $g$-torus $\CE$ is an abelian variety if and only if
it is bi-holomorphic with a torus $\C^g/\Lambda$ with  $\Lambda=\tau \Z^g+ D\Z^g$,
where $D$ is a diagonal matrix
$D=Diag(d_1,\dotsc,d_g)=\begin{pmatrix} d_1 &  & \text{\hspace{-5pt}{\Large $0$}}\\
& \ddots \\
 \text{\Large $0$} & & d_g
\end{pmatrix}$
with positive integers $d_1|d_2|\dotsb |d_g$ and $\tau$ is a complex $g\times g$
symmetric matrix with a positive definite imaginary part.
\end{thm}

\begin{defn}\
\begin{enumerate}
\item The set of matrices
\[
 \CH_g =\{ \tau \in \Sym(g,\C) \colon Im(\tau) \text{ is positive definite }\}
\] is called
 {\em the Siegel
  upper half-space of degree $g$}. \vlabel{note:siegel}
\item  An integer diagonal matrix $D=Diag(d_1,\dotsc,d_g)$ with
positive $d_1,\dotsc,d_g$ satisfying $d_1|d_2|\dotsc|d_g$ is
called {\em a polarization type}. When $D=I_g$ we say that the
polarization is {\em principal}.
 \item We say that a
$g$-torus $\CE$ admits a polarization of the type $D=Diag(d_1,\dotsc,d_g)$ if $\CE$
is bi-holomorphic with  a torus $\C^g/\Lambda$ with $\Lambda=\tau\Z^g+D\Z^g$ for
some $\tau\in \CH_g$.
\end{enumerate}
\end{defn}

Thus a torus $\CE$ is an abelian variety if and only is it admits a polarization.

\begin{rem}\vlabel{manyD} It is not hard to see that every abelian variety $\C^g/\Lambda$ admits infinitely
  many polarization. For example, since for every positive integer $k$
  abelian varieties $\C^g/\bigl(\tau\Z^g+D\Z^g\bigr)$ and
  $\C^g/\bigl(k\tau\Z^g+kD\Z^g\bigr)$ are isomorphic (via the map
  $z\mapsto kz$), every abelian variety admitting a polarization of type $D$ also admits a
polarization of type $kD$ for any positive $k\in \N$.
\end{rem}

\begin{ntn}\vlabel{note:polar} For a polarization type $D$ and $\tau\in \CH_g$ we will denote
by $\Lambda^D_\tau<\C^g$ the lattice
\[ \Lambda^D_\tau=\tau\Z^g+D\Z^g,\]
and by $\pav D \tau$ the
 abelian variety \
\[\pav D \tau =\C^g/(\tau\Z^g+D\Z^g).\]
\end{ntn}

\subsection{Action of $Sp(2g,\R)$ and isomorphisms of polarized abelian varieties}
For $K=\R,\Q$, or $\Z$ we will denote by $Sp(2g,K)$ the corresponding
symplectic group\vlabel{note:sp}
\[ Sp(2g,K) =\left\{ M\in Gl(2g,K) \colon
M \begin{pmatrix} 0 & I_g \\ -I_g & 0 \end{pmatrix} \tra M =
\begin{pmatrix} 0 & I_g \\ -I_g & 0 \end{pmatrix}\right\}. \]
It is  a subgroup of $Gl(2g,K)$ closed under transposition.

The group $Sp(2g,\R)$ acts (on the left)  on the Siegel upper
half-space $\CH_g$ via
\[ \begin{pmatrix} \alpha & \beta \\
  \gamma & \delta \end{pmatrix} \cdot \tau = (\alpha\tau
+\beta)(\gamma \tau + \delta)^{-1}.\]

For a polarization type $D$, let $G_D$ be the following subgroup
of $Sp(2g,\Q)$:\vlabel{note:Gd}
\[
G_D= \left\{ M\in Sp(2g,\Q) \colon
\begin{pmatrix} 1 & 0 \\ 0 & D \end{pmatrix}^{-1}
M
\begin{pmatrix} 1 & 0 \\ 0 & D \end{pmatrix} \in M_{2g}(\Z)\right\}. \]

Notice that $G_{I_g}=Sp(2g,\Z)$.

Without going into details of polarizations we just state the
following two facts that follow from \cite[Proposition 8.1.3 and
Remark 8.1.4]{BL}.

\begin{fact}\vlabel{mainfact} Let $D$ be a polarization type, $A\in GL(g,\C)$, and
$\tau,\tau'\in \CH_g$. The map $z\mapsto Az$ induces an
isomorphism (as polarized abelian varieties) between  the
polarized abelian varieties $\pav{D}{\tau'}$ and $\pav{D}{\tau}$
if and only if there is $M=\begin{pmatrix} \alpha & \beta \\
\gamma & \delta
\end{pmatrix}\in G_D$ such that $\tau' = M\cdot \tau$ and $A=\tra{(\gamma\tau+\delta)}$.
\end{fact}

\begin{rem}\vlabel{remain} For $A\in GL(g,\C)$ the matrix $M$ in
the previous fact is the unique matrix $M\in Gl(2g,\R)$ such that
$A(\tau,,I_g)=(\tau,I_g)\tra M$, namely  the following diagram is
commutative.
\[%
\begin{diagram}\node{\R^{2g}}
\arrow[5]{e,t}{x\mapsto (\tau'\, ,\, I_g)x }
\arrow[2]{s,l}{x\mapsto \tra M x} \node[5]{\C^g}
\arrow[2]{s,r}{z\mapsto Az}\\[2]
\node{\R^{2g}} \arrow[5]{e,b}{x\mapsto (\tau\,, \, I_g)x }
\node[5]{\C^g}
\end{diagram}
\]
\end{rem}

\subsection{The Siegel fundamental set $\fF_g$ for the action of $Sp(2g,\Z)$}
We fix a positive $g\in \N$  and let $n=g(g+1)/2$. We will identify
 $\Sym(g,\R)$  with $\R^n$,
and $\Sym(g,\C)$  with $\C^n$.

\subsubsection{Minkowski reduced matrices}

\begin{defn}\vlabel{def:minred}
 A real symmetric matrix $\beta=(\beta_{i,j})\in \Sym(g,\R)$
 is called {\em Minkowski reduced}  if it satisfies the following conditions.
\begin{enumerate}[M(I)]
\item $\beta_{1,1}>0$; \item $\beta[a]\geq \beta_{k,\, k}$, for
all $k=1,\dotsc g$, and   for all $\tra a=(a_1,\ldots, a_g)\in
\mathbb Z^g$, with $\gcd(a_k,\ldots, a_g)=1$. \item $\beta_{k,\,
k+1}\geq 0$ for all $k=1,\ldots, g-1$.
\end{enumerate}

We denote by $\Mi_g$  the set of all Minkowski-reduced  $g\times g$ real matrices.
\end{defn}

\begin{rem}\vlabel{rem32}
  Note that the definition above differs with the definition in
  \cite[p.191]{Igusa}, or in \cite[p.128]{TerrasII} in that condition
  M(I) here is replaced there by the assumption that the matrix $\beta$
  is positive definite. However, if a matrix $\beta$ is positive
  definite then $\beta _{1,1}>0$ and conversely, Siegel shows in
  \cite[Theorem 3]{oldSiegel} that every matrix satisfying M(I), M(II)
  is positive definite, so the definitions are equivalent.
It follows  that the set $\Mi_g$ is closed in the set of all real
symmetric positive-definite $g\times g$-matrices.
\end{rem}

Here are some basic facts about  Minkowski reduced matrices (see
references below).

\begin{fact}\vlabel{M-reduced} \begin{enumerate}[(a)]
\item  There are positive real constants $c,c'$ such that for
every  real matrix $\beta\in \mathcal \Mi_g$, and for every $x
\in\R^g$,
$$c \sum_{i=1}^g\beta_{i,\,i}x_i^2\leq \beta[x] \leq c' \sum_{i=1}^g\beta_{i,\,i}x_i^2$$

\item For every matrix $\beta\in \mathcal \Mi_g$ and $1\leq i,j\leq
  g$, we have $|\beta_{i,\,j}|\leq \frac{1}{2}\beta_{i,\,i}.$

\item  For every matrix $\beta\in \mathcal \Mi_g$, we have $0 <
\beta_{1,\,1}\leq \beta_{2,\,2}\leq \cdots \leq \beta_{g,\,g}.$

\item There are finitely many inequalities in M(II), which together
with M(I) and M(III) define $\Mi_g$. In particular, the set $\Mi_g$ is semi-algebraic.

\item  The dimension of $\Mi_g$ is $n$.

\end{enumerate}
\end{fact}
 For (a), (b), and  and (c) above, see \cite[p. 132, Proposition 1(c,a)]{TerrasII}.
Clause (d) is exactly \cite[p. 130, Theorem 1(2)]{TerrasII}.
For $(e)$, it follows from \cite[Theorem 6]{oldSiegel} that the
$\Mi_g$ contains a non-empty open subset of positive definite $g\times
g$ matrices. It is easy to see that the set of positive definite matrices is open
in $\Sym(g,\R)$, hence the dimension of $\Mi_g$ is $n$.

\subsubsection{The Siegel fundamental set}
For the following definition see \cite[p. 194]{Igusa}.

\begin{defn} The Siegel fundamental set  $\fF_g$ is  the collection  of all
symmetric $\tau\in Sym(g,\C)$ such that
\begin{enumerate}[(i)]
\item $Im(\tau)$ is in $\Mi_g$. \item $|Re(\tau_{i,j})|\leq 1/2$
for all $1\leq i,j\leq g$. \item $det (Im(\sigma \tau))\leq
det(Im(\tau))$ for all $\sigma\in Sp(2g,\Z)$.
\end{enumerate}
\end{defn}

As is pointed out in  \cite[p.194]{Igusa}, clause (iii) is equivalent to:
 \begin{equation}\elabel{equation1} \mbox{For all }\sigma=
\begin{pmatrix}
  \alpha & \beta \\
  \gamma & \delta \\
\end{pmatrix}
\in Sp(2g,\Z),
 \mbox{ we have } |det(\gamma\tau+\delta)|\geq 1. \end{equation}

Notice that, since every Minkowski-reduced matrix is positive-definite, we have
$\fF_g\subseteq \CH_g$.

Here are the main properties that we are going to use.

\begin{fact}\vlabel{Igusa15}\leavevmode
\begin{enumerate}[(a)]
\item $\fF_g$ is a closed subset of $\C^n$ (and not only of
$\CH_g$). \item If $\tau=\alpha+i\beta\in \fF_g$ then
$\beta_{1,\,1}\geq \frac{\sqrt{3}}{2}.$ \item There are only
finitely many inequalities in (III), which together with (I) and
(II) define $\fF_g$. In particular $\fF_g$ is a semi-algebraic
subset of $\CH_g$.
\end{enumerate}
\end{fact}
For  item (a) see \cite[Lemma 16, p. 196]{Igusa}. For
(b) see \cite[Lemma 15, p. 195]{Igusa}, and for (c) we refer to
\cite[Theorem 1, p.35]{klingen}.

The following fact explains the significance of the set $\fF_g$.
\begin{fact}\vlabel{Igusa15.1} The set
$\fF_g\sub \C^n$ is the Siegel fundamental set for the action of $Sp(2g,\Z)$ on
$\mathcal H_g$, namely
\begin{enumerate}[(a)]
\item
 $Sp(2g,\Z)\cdot \fF_g=\CH_g$ and the there are only finitely
many elements $\sigma\in Sp(2g,\Z)$ such that $\fF_g\cap \sigma\cdot \fF_g$ is nonempty.
\item For every compact $X\subseteq \CH_g$ the set $\{ \sigma\in
  Sp(2g,\Z)\colon \sigma\cdot \fF_g\cap X \not=\varnothing\}$ is finite.
\end{enumerate}
\end{fact}
See  \cite[Theorem 2, p.34]{klingen} for the proof of the above
fact.
\medskip

We now describe the fundamental set $\fF_g(G)$, for  $G<Sp(2g,\Z)$
a subgroup of finite index. Let $\gamma_1,\dots, \gamma_k\in
Sp(2g,\Z)$ be representatives of all left co-sets of $G$ in
$Sp(2g,\Z)$ (so $Sp(2g,\Z)=\cup_{i=1}^k G\gamma_i$). We call a set
of the form $\cup_{i=1}^k \gamma_i\cdot \fF_g$ {\em a Siegel
fundamental set for $G$},
 and suppressing the dependence on representatives, we will denote this set by $\fF_g(G)$.
So  \vlabel{note:fFD}
\begin{equation} \label{fFD} \fF_g(G)=\cup_{i=1}^k \gamma_i\cdot \fF_g. \end{equation}

Obviously, $\fF_g(G)$ is a semi-algebraic subset of $\CH_g$, and
since $\fF_g$ is a fundamental set for the action of $Sp(2g,\Z)$,
the set $\fF_g(G)$  is a fundamental set for the action of $G$ on
$\CH_g$ and we have:
\begin{fact}\vlabel{fund-set}\begin{enumerate}[(i)]
\item  $\fF_g(G)$ is a closed
subset of $\CH_g$;
\item  $G\cdot \fF_g(G)=\CH_g$;
\item there
are only finitely many  $\sigma\in G$ such that
$\fF_g(G)\cap \sigma\cdot \fF_g(G)$ is nonempty.
\item For every compact $X\subseteq \CH_g$ the set $\{ \sigma\in
  G\colon \sigma\cdot \fF_g(G)\cap X \not=\varnothing\}$ is finite.
\end{enumerate}

\end{fact}

\section{The classical theta functions and embeddings of abelian varieties}

\subsection{The theta functions}
As pointed out earlier,  we identify the set $Sym(g,\C)$ with
$\C^n$ (for $n=g(g+1)/2$), and view   $\CH_g$ as a subset of
$\C^n$.

For $(z,\tau)\in \C^g\times \CH_g$ let  \vlabel{note:theta}
\[ \vartheta(z,\tau) = \sum_{n\in \Z^g} \eexp{ \pi i(\tra n\tau n +2\,\tra n z)}.   \]
It is known (see for example \cite[p. 118]{Mumford}) that the
above series is convergent and $\vartheta$ is holomorphic on
$\C^g\times \CH_g$. It is also immediate from the definition that
$\vartheta(z,\tau)$ is $\Z^g$-periodic in $z$ and
$(2\Z)^n$-periodic in $\tau$.

\begin{defn}  For $a,b\in
\R^g$, the associated {\em Riemann Theta function} is the function $\teta{a}{b}
: \C^g\times \CH_g\to \mathbb C$ defined by:
\[ \teta a b (z,\tau)=\eexp{\pi i(\tra a \tau  a + 2\, \tra a(z+b))}\vartheta(z+\tau a +b,\tau).
\]
\end{defn}

\subsection{Embeddings of abelian varieties}
We fix a polarization type $D=Diag(d_1,\dotsc,d_g)$. We also fix a
set of representatives $\{ c_0,\dotsc,c_N\}$ of the cosets of
$\Z^g$ in the group $D^{-1}\Z^g$.

\begin{rem} We could take $\{c_0,\dotsc,c_N \}$ to be the set of all vectors $c\in D^{-1}\Z^g$
whose components lie in the interval $[0,1)$. In particular we have
$N=d_1\cdot d_2\cdot \dotsb d_g-1$.
\end{rem}

The following fact is a consequence of the classical Lefschetz
Theorem (see \cite[Theorem 1.3 p.128]{Mumford})
\begin{fact}[Lefschetz Theorem]\vlabel{lef}
\leavevmode
\begin{enumerate}[(a)]
\item Assume  $d_1\geq 2$. \\The functions
$\left\{ \teta {c_0} 0(z,\tau), \dotsc, \teta {c_N} 0(z,\tau)\right\}$ have no zero in common, and
\[ \varphi^D (z,\tau)= \left(\teta {c_0} 0(z,\tau)\,
:\, \teta {c_1} 0(z,\tau) :\,
 \dotsc\, :\,\teta {c_N} 0(z,\tau) \right) \]
is a well-defined holomorphic map from $\C^g\times \CH_g$ into $\PP^N(\C)$. \\
For each $\tau\in \CH_g$ the map $\varphi^D_\tau\colon z\mapsto
\varphi^D(z,\tau)$ is $\Lambda^D_\tau$-periodic, hence induces a
holomorphic map from the abelian variety $\CE^D_{\tau}$ into
$\PP^N(\C)$. \item Assume $d_1\geq 3$. \\For each $\tau\in \CH_g$
the map $\varphi^D_\tau$ is an immersion on $\C^g$ that induces an
analytic
  embedding of the abelian variety $\pav D \tau$ into the
  projective space $\PP^N(\C)$, whose image is an algebraic variety.
\end{enumerate}
\end{fact}

From now on, in the case $d_1\geq 2$, we will denote by $\varphi^D(z,\tau)$ the map \vlabel{note:vd}
\[ \varphi^D (z,\tau)= \left(\teta {c_0} 0(z,\tau)\,
:\, \teta {c_1} 0(z,\tau) :\,
 \dotsc\, :\,\teta {c_N} 0(z,\tau) \right), \]
and for $\tau\in \CH_g$, we will denote by $\varphi^D_\tau(z)$ the
map \vlabel{note:vdt} $\varphi^D_\tau\colon z\mapsto
\varphi^D(z,\tau)$ from $\C^g$ into $\PP^N(\C)$.

\section{Definability of holomorphic $\Z^n$-invariant maps}
The main result is in this section is a general theorem (see
Theorem \ref{defpseta} below) about definability in the structure
$\R_{an,exp}$ of certain periodic holomorphic functions on
``truncated'' tube domains. The theorem will then be applied to
prove the definability of the theta functions on a restricted
domain. We first review the basics of polyhedral cones and tube
domains.

\subsection{Polyhedral cones and Tube domains}
\begin{ntn}
 For vectors $v,u \in \R^m$ we write $v > u$ ($v\geq u$)
if $v_i>u_i$ ($v_i\geq u_i$) for all $i=1,\dotsc,m$.
\end{ntn}

Let $\FF$ be the field of complex or real number, and $L\colon \FF^m\to \FF^n$ be a
linear map. We say that $L$ is {\em integral} if the standard matrix of $L$ (i.e.
the matrix representing $L$ with respect to the standard bases) has integer
coefficients. Clearly,  $L$ is integral if and only if it maps $\Z^m$ into $\Z^n$.

An nonempty subset $C\subseteq \R^n$ is called {\em a
polyhedral
  cone} if $C=\{ x\in \R^n \colon Ax\geq 0\}$ for some $m\times n$-matrix
$A$. The cone $C$ is called {\em integral} if one can choose $A$ in $M_{m\times
n}(\Z)$. The dimension $\dim(C)$ of a polyhedral cone $C$ is the dimension of its
linear span.  We call $(\R^{\geq 0})^n$ {\em the standard polyhedral cone} in $\R^n$
and denote it by $\spc n$. \vlabel{note:spc}

By the Minkowski-Weyl theorem (see for example \cite[Theorem 1.3]{Ziegler}), for every  polyhedral cone $C\subseteq \R^n$, there
are $v_1,\ldots, v_k\in C$, such that
\[ C=\left\{\sum\lambda_j v_j:\lambda_j\geq 0 \right\}.\]
Moreover, since Fourier-Motzkin elimination used in the proof of \cite[Theorem 1.3]{Ziegler} works
over $\Q$,
 if $C$ is integral then $v_1,\ldots, v_k$ can be
chosen in $\mathbb Z^n$. Notice, that in this case
$\dim(C)$ is  the linear dimension of $\{v_1,\ldots, v_k\}$.

Stated differently,

\begin{fact}\vlabel{Minkowski}
  If  $C\subseteq \R^n$ is a  polyhedral cone of dimension $k$ then
  there is a linear map $L\colon \R^m\to \R^n$ of rank $k$ such that
  $C$ equals the image of $\spc m$ under $L$.
Moreover, if $C$ is integral then $L$ can be chosen to be integral as well.
\end{fact}

\begin{rem}\vlabel{conint}
In the above setting, the relative interior of $C$ (i.e. the interior of $C$ in its
linear span) is $L((\R^{>0})^m)$ (see for example \cite[Proposition 2.1.12]{Fundca}
).
\end{rem}

\begin{defn}
For  $S\subseteq\R^n$ the following complex subset  is  called
{\em the tube domain associated to $S$}:
\[ \CT(S)=\{x+iy\in
\C^n:   y\in S\}.\] Also, for a positive real $R>0$ we will denote by $\CT_{\leq
R}(S)$ {\em the truncated tube domain}
\begin{multline*}
\CT_{\leq R}(S)=
\{x+iy\in\C^n: \bigwedge_{i=1}^n |x_i|\leq R,\,  y\in S\}=\\
\{ z\in \CT(S)\colon \bigwedge_{i=1}^n |Re(z_i)|\leq R \}.
\end{multline*}
\end{defn}

\begin{ntn}
For a linear map   $L\colon \R^k\to \R^n$  we
 denote by $L_\C$  the complexification of $L$, i.e. the $\C$-linear  map $L_\C\colon
\C^k\to \C^n$ defined by $L_\C(x+iy)= L(x)+iL(y)$ for $x,y\in \R^k$.
\end{ntn}

\begin{claim}\vlabel{Lonto} Let $S_1\sub \R^m, S_2\sub \R^n$ be sets
and   $L:\R^m\to \R^n$  a linear map with $L(S_1)=S_2$. If $S_2$
has non-empty interior then $L_\C$ maps $\CT(S_1)$ onto
$\CT(S_2)$.
\end{claim}
\proof Since the image of $L$ contains an open subset of $\R^n$ its rank equals $n$,
and therefore $L(\R^m)=\R^n$. The result now follows from  the definition of a tube
domain. \qed

\subsection{Definability of holomorphic $\Z^n$-invariant maps}
In this section we prove the following general theorem about
definability in the structure $\R_{an,exp}$ of certain periodic
holomorphic functions on truncated tube domains.

\begin{ntn}\vlabel{cvfd} Let $C\subseteq \R^n$ be a subset.
\begin{enumerate}
\item For a  vector $v\in \R^n$, we denote by $C\la v\ra$ the translate of $C$ by
the vector $-v$, namely
\[ C\la v\ra = C-v =\{ x\in \R^n \colon x+v\in C \}, \]
\item For a function $f\colon \R^n\to \R$  and a real number $d\in \R$, we use the
following notations:
\[ C^{f>d} = \{ x\in C   \colon f(x)> d \},\quad C^{f\geq d} = \{ x\in C   \colon f(x)\geq  d \}.\]
\end{enumerate}
In particular $C\la v\ra^{f>d}$ denotes the set
\[ C\la v \ra^{f>d} = \{ x\in \R^n\colon x+v\in C, f(x)> d \}.\]
\end{ntn}

\begin{thm}\vlabel{defpseta}
Let $C\sub \R^n$ be an integral  polyhedral cone  with
$\dim(C)=n$,  let $\ell:\R^n\to \R$ be  an integral linear
function  which is positive on $Int(C)$,  let $d_0>0 $ be a
positive real number, and $v\in Int(C)$.

Let $U\supseteq \CT\bigl(C\la v\ra^{\ell> d_0}\bigr)$ be an open subset of $\C^n$,
and $\theta\colon U\to \C$ a holomorphic, $\Z^n$-periodic function bounded on
$\CT\bigl(C\la v\ra^{\ell> d_0}\bigr)$.

Then, for any real $d>d_0$,  $u\in C$ with $u<v$, and for any real
$R>0$,
  the restriction
  of $\theta$ to the closed set
 $\CT_{\leq R}\bigl(C\la u\ra ^{\ell \geq d}\bigr)$
is definable in the
  structure $\R_{an,exp}$.

Moreover,  there is a definable open $V\sub
  \C^n$ with  $V\supseteq \CT_{\leq R}\bigl(C\la u\ra^{\ell\geq d}\bigr)$,
such that the restriction  $\theta\rest V$ is
  definable in $\R_{an,exp}$.
\end{thm}
\proof We fix $d>d_0>0$ and $u<v$ in $C$ as in the theorem. We
first consider a  special case.

\medskip
\noindent{\bf The special case: $C=\spc n$} \\
 Notice that because of the
$\Z^n$-periodicity of $\theta$, it is sufficient to prove the result for $R=1/2$.

We use $\e(z)$ to denote the complex exponential function $\e(z)=e^{2\pi i z}$
\vlabel{note:be}
and we
use $\e_n:\C^n \to \C^n$ to denote the map which is $\e(z_i)$ in each
coordinate. The map $\e_n$ is a surjective group homomorphism from $\la \C^n,+\ra$
to $\la (\C^*)^n,\cdot\ra $ which is  locally a bi-holomorphism, with kernel
$\Z^n$.

Since $\theta$ is $\Z$-periodic, it factors through the map $\e_n$, i.e. there is a
holomorphic map $\theta^*\colon \e_n(U)\to \C$
 such that the following diagram commutes.
\[
\begin{diagram}
\node[2]{U}\arrow{sw,t}{\theta}\arrow{s,l}{\e_n}\\
 \node[1]{\C} \node{\e_n(U)}\arrow{w,t}{\theta^*}
\end{diagram}
\]

Since $d>d_0$, and $u<v$  we have
$C\la u\ra^{\ell\geq d} \subseteq C\la v\ra^{\ell>d_0}$, hence
$\CT\bigl(C\la u\ra^{\ell \geq d}\bigr) \sub U$.

To simplify notations, let
$T_v^{>d_0}=\CT\bigl(C\la v\ra^{\ell> d_0}\bigr)$,
$T_u^{\geq d}=\CT_{\leq 1/2}\bigl(C\la u\ra^{\ell\geq d}\bigr)$,
$\wt T_v^{>d_0} =\e_n(T_v^{>d_0})$, and $\wt T_u^{\geq d}=\e_n(T_u^{\geq d})$.

We have a commutative diagram where all vertical arrows are surjective maps.
\[
\begin{diagram}
\node[2]{T_v^{>d_0}}\arrow{sw,t}{\theta}\arrow{s,l}{\e_n}\node{T_u^{\geq
d}}\arrow{w,t,L}{\id}
\arrow{s,l}{\e_n}\\
 \node[1]{\C} \node{\wt T_v^{>d_0}}\arrow{w,t}{\theta^*}\node{\wt T_u^{\geq d}}
\arrow{w,t,L}{\id}
\end{diagram}
\]

Note that the restriction of $\e(z)$ to the set $\{ z\in \C \colon |Re(z)|\leq
1/2\}$ is definable in $\Rae$, since it only requires the functions $\exp(x)\colon
\R\to \R$, $\sin(x)\rest[-\pi,\pi]$, and $\cos(x)\rest[-\pi,\pi].$ Hence, the
restriction of $\e_n$ to $T_u^{\geq d} $ is definable in $\Rae$ and it is sufficient
to show that the restriction of $\theta^*$ to $\wt T_u^{\geq d}$
is definable in $\R_{an}$. \\

Let $v_1,\dotsc,v_n\in \R$ be the components of $v$, i.e.
$v=(v_1,\dotsc,v_n)$, and $l_1,\dotsc,l_n\in \Z$ be the components
of $\ell$, i.e. $\ell(x)=l_1x_1+\dotsb+l_nx_n$. Notice that since,
by assumptions, $\ell$ is positive on $Int(\spc n)$, all $l_i$ are
positive integers. Let  $u_1,\dotsc u_n\in \R$ be the components
of $u$. We have $0\leq u_i<v_i$ for $i=1,\dotsc n$.

Recall that
\begin{gather*}
T_v^{>d_0}=\{ x+iy\in \C^n \colon   \bigwedge_{i=1}^n  y_i \geq  -v_i, \quad
l_1 y_1+\dotsb +l_n y_n > d_0 \} \text{ and }\\
 T_u^{\geq d} =\{ x+iy\in \C^n \colon \bigwedge_{i=1}^n |x_i| \leq 1/2,\,  \bigwedge_{i=1}^n  y_i \geq
 -u_i\, ,\,
l_1 y_1+\dotsb +l_n y_n \geq d \}.
\end{gather*}

It is not hard to compute the images of $T_v^{>d_0}$ and $T_v^{>d_0}$  under $\e_n$ and obtain
\begin{gather*}
 \wt T_v^{>d_0}= \{ q\in \C^n \colon
\bigwedge_{i=1}^k 0<|q_i|\leq e^{2\pi v_i},\, |q_1|^{\l_1}\cdots
|q_n|^{\l_n} <  e^{-2\pi d_0}\},\\
 \wt T_u^{\geq d}= \{ q\in \C^n \colon
\bigwedge_{i=1}^k 0<|q_i|\leq e^{2\pi u_i},\, |q_1|^{\l_1}\cdots
|q_n|^{\l_n} \leq  e^{-2\pi d}\}.
\end{gather*}

Let
\[ O=\{ z\in \C^n\colon \bigwedge_{i=1}^n |q_i|< e^{2\pi v_i},\,
|q_1|^{l_1}\dotsb |q_n|^{l_n} <e^{-2\pi d_0} \}\]
and
\[Z=\{ q\in \C^n\colon \prod_{i=1}^{n}q_i=0\}.\]
Obviously  $O$ is an open subset of $\C^n$, and $O \setminus Z
\subseteq \wt T_v^{>d_0}$. Since $\theta$ is bounded on
$T_v^{>d_0}$, the function $\theta^*$ is bounded on $O\setminus
Z$, and therefore,   by Riemann's removable of singularities
theorem for complex functions of several variables, $\theta^*$ has
holomorphic extension $\Theta^*\colon O\to \C$.

Consider the set
\[ W= \{ q\in \C^n \colon
\bigwedge_{i=1}^k |q_i|\leq e^{2\pi u_i},\, |q_1|^{\l_1}\cdots
|q_n|^{\l_n} \leq  e^{-2\pi d}\}.
\]
Obviously $W$  is a  closed subset of $\C^n$, and, since all $l_i$ are positive, it
is also bounded, hence compact. Since $u_i<v_i$, and $d>d_0$, we have $W\subseteq
O$. Clearly, $W$ is a semi-algebraic set, therefore the restriction of $\Theta^*$ to
$W$ is definable in $\Ran$.

The set $\wt T_u^{\geq d}=W\setminus Z$ is semialgebraic. Hence
the restriction of $\Theta^*$ to $\wt T_u^{\geq d}$ is definable
in $\Ran$.  This finishes the proof of the special case.
\medskip

\noindent{\bf General Case}

Let $L\colon \R^k \to \R^n$ be an integral linear map, as in Claim \ref{Minkowski},
sending $\spc k$ onto the cone $C$. By Remark 5.3 there is $v'\in Int(\spc k)$ with $L(v')=v$.
Since $L$ is a linear map it maps $\spc k\la v'\ra$ onto $C\la v\ra$, and, by Claim \ref{Lonto},
 $L_\C$ maps $\CT(\spc k)$ onto $\CT(C)$, and
$\CT(\spc k\la v'\ra)$ onto $\CT(C\la v\ra)$.

\medskip
Let $\ell'=\ell\co L$, $\theta' = \theta\co L_\C$, and $U'=L_\C^{-1}(U)$.

Using Claim \ref{Lonto} again we obtain that $L_\C$ maps $\CT(\spc k\la
v'\ra^{\ell'>d_0})$ onto $\CT(C\la v\ra^{\ell >d_0})$, hence $U'$ contains $\CT(\spc
k\la v'\ra^{\ell'>d_0})$ and $\theta'$ is bounded on $\CT(\spc k\la
v'\ra^{\ell'>d_0})$.

Because $\theta'$ is $\Z^n$-periodic, we can now use the special case and obtain
that the restriction of $\theta'$ to the set $\CT_{\leq R'}(\spc k\la
u'\ra^{\ell'\geq d})$ is definable in $\Rae$. It follows then that the restriction
of $\theta$ to $\CT_{\leq R}(C\la u\ra^{\ell \geq d})$ is definable as well.

This finishes the proof of the main part of the theorem.\\

For the ``moreover'' part of the statement, we can take $V$ to be
the interior of the set $\CT_{\leq R'}(C\la u'\ra^{\ell \geq d'})$
for any $R'> R$, $d_0<d'<d$, and $u'\in C$ with
$u<u'<v$. \\

This  concludes the proof of of Theorem \ref{defpseta}. \qed

\section{Definability of  theta functions}

Our ultimate goal in this section is to show that the restrictions of
Riemann Theta functions $\teta{a}{b}(z,\tau)$ to an appropriate
sub-domain of $\C^g\times \CH_g$ is definable in $\Rae$.  Towards this
goal we need to establish the assumptions of Theorem \ref{defpseta},
and mainly establish the boundedness of some variations of
$\vartheta(z,\tau)$.  We use ideas from similar boundedness proofs in
\cite{Mumford}.

\medskip

We fix a positive integer $g$ and let $n=g(g+1)/2$. We identify
the set of all real symmetric $g\times g$-matrices $Sym(g,\R)$
with $\R^n$. We will also identify $Sym(g,\C)$ with $\C^n$ and
with $Sym(g,\R)+iSym(g,\R)$.  In particular for a set $S\subseteq
Sym(g,\R)$ we view the corresponding cone $\CT(S)$ as a subset of
$\Sym(g,\C)$.

\begin{defn}\vlabel{note:miclosed} Let $\Mib_g\subseteq Sym(g,\R)$ be the integral polyhedral cone given by
$\beta_{1,1} \geq 0$, $\beta_{k,k+1}\geq0, k=1,\dotsc,g-1$, and
the finitely many inequalities needed for the condition M(III) in
Definition \ref{def:minred} (See Fact \ref{M-reduced}(d)).
\end{defn}

It is immediate from the definition that
\[\Mi_g=\Mib_g \cap \{ \beta\in Sym(g,\R) \colon
\beta_{1,1}>0\},\]
and the interior of $\Mib_g$ is contained in $\Mi_g$. In particular the dimension of $\Mib_g$ is $n$.

\begin{ntn}\vlabel{note:notec}
For a real $m>0$, we denote by $C_m$  the following integral polyhedral cone
in $\R^g\times \R^n$:
\[
C_m=\left\{(y,\beta)\in \R^g\times Sym(g,\R) \colon
\beta\in \Mib_g \, , \, \bigwedge_{i=1}^g |y_i|\leq m\beta_{i,i}\right\}.
\]
\end{ntn}
Let $\ell_{1,1}\colon \R^g\times \R^n\to $ be the integral linear function
$\ell_{1,1}\colon (y,\beta) \mapsto \beta_{1,1}.$

Notice that for every $d>0$ we have $\CT(C_m^{\ell_{1,1}\geq
d})\subseteq \C^g\times \CH_g$.

\medskip
Before proving definability of $\teta{a}{b}$ we need to establish
the boundedness of an auxiliary function.

\begin{prop}\vlabel{prop1} Let $m,d>0$ be real numbers. Then there is a positive $k\in
\mathbb N$ and $\beta^*\in Int(\Mib_g)$  such that for
$\beta_0^*=(0,\beta^*)$ the set $\CT(C_m\la
\beta^*_0\ra^{\ell_{1,1}>d})$ is contained in $\C^g\times \CH_g$
and
 the function
\[ \theta (z,\tau)=\eexp{2\pi i k \tau_{g,g}}\vartheta(z,2\tau) \]
is holomorphic and bounded on $\CT(C_m\la \beta^*_0\ra^{\ell_{1,1}>d})$.

\end{prop}

\proof
Recall that
\begin{multline*}
\CT(C_m\la \beta^*_0\ra^{\ell_{1,1}>d}) =\\
\{(z,\tau) \in \C^g\times\C^n \colon (Im(z),Im(\beta)+\beta^*)\in C_m,\,
Im(\tau_{1,1}) > d\}.
\end{multline*}

Let $c,c'$ be real constants as in Fact \ref{M-reduced}(a).
Namely, for every  real matrix $\beta\in \mathcal \Mi_g$, and for
every $x \in\R^g$,
$$c \sum_{i=1}^g\beta_{i,\,i}x_i^2\leq \beta[x] \leq c' \sum_{i=1}^g\beta_{i,\,i}x_i^2.$$

We take $k\in \Z$ and $\beta^*\in Int(\Mib_g)$ satisfying
\[ k> \frac{m^2g}{2c} \text{\ \  and\ \  } \beta^*_{g,g}< \frac{cd}{2c'}.\]
Notice that  since $\Mib_g$ is a cone of dimension $n$  such $\beta^*$ exists.

Since $\beta^*\in Int(\Mib_g)$, we have $\beta^*\in \Mi_g$, hence by
Fact \ref{M-reduced}(c)
\[ 0<\beta^*_{1,1}\leq \dotsc \leq \beta^*_{i,i}\leq \dotsc\leq  \beta^*_{g,g}< \frac{cd}{2c'}. \]

Let $(z,\tau)\in \CT(C_m\la \beta^*_0\ra^{\ell_{1,1}>d})$ with $z=x+iy$. For
$\beta=Im(\tau)+\beta^*$ we have $(y,\beta)\in C_m$.

Since $(x+iy,\tau)\in \CT(C_m\la\beta^*_0\ra^{\ell_{1,1}>d})$, we have $\beta_{1,1}-\beta^*_{1,1}>d$, hence  $\beta_{1,1}>d$,
$\beta\in \Mi_g$, and, by Fact \ref{M-reduced}(c),
\[ d< \beta_{1,1}\leq \dotsc \leq \beta_{i,i}\leq \dotsc\leq  \beta_{g,g}. \]

For any $v\in \R^g$ we have
\begin{equation}\elabel{est1}
\begin{split}
  Im(\tau)[v] &=
  (\beta-\beta^*)[v]=\beta[v]-\beta^*[v] \\
&\geq c\sum \beta_{i,i}v_i^2 -c' \sum \beta_{i,i}^*v_i^2\\
 & =\sum(c\beta_{i,i}-c'\beta^*_{i,i}) v_i^2 \geq
  \sum\left(c\beta_{i,i} -c'\frac{cd}{2c'}\right)v_i^2 \\
  &=\sum\left(c(\beta_{i,i} -\frac{d}{2})\right)v_i^2
\geq
   \frac{c}{2}\sum \beta_{i,i}v_i^2.
\end{split}
\end{equation}
Thus $Im(\tau)$ is a positive definite symmetric matrix, hence $\tau\in
\CH_g$ and $\CT(C_m\la \beta^*_0\ra^{\ell_{1,1}>d})$ is contained in $\C^g\times \CH_g$.

Since $\vartheta$ is holomorphic on $\C^g\times \CH_g$, we
obtain that $\theta$ is also holomorphic on $\CT(C_m\la \beta^*_0\ra^{\ell_{1,1}>d})$.

To  show the boundness of $\theta$
on $\CT(C_m\la \beta^*_0\ra^{\ell_{1,1}>d})$ we
use ideas form \cite[p. 118 Proposition 1.1.]{Mumford}.
\begin{equation}\elabel{est2}
\begin{split}
|\theta (z,\tau)|
&=\bigl|\eexp{2\pi i k \tau_{g,g}}\bigr|\cdot
\bigl|\sum_{n\in \Z^g} \eexp{ \pi i(2\tau[n] +2\,\tra n z)}\bigr| \\
&\leq
|\eexp{2\pi i k \tau_{g,g}}|\cdot
\sum_{n\in \Z^g}\bigl| \eexp{ \pi i(2\tau[n] +2\,\tra n z)}\bigr|
\\
&= \eexp{-2\pi k Im(\tau_{g,g})} \cdot
\sum_{n\in \Z^g} \eexp{ -\pi(\,  2Im(\tau)[n] -\,2\tra n y\,)}\\
&= \sum_{n\in \Z^g} \eexp{ -\pi\bigl(\,  2Im(\tau)[n] -\,2\tra n
y\,+2k Im(\tau_{g,g})\,\bigr)}.
\end{split}
\end{equation}

Because $(y,\beta)\in C_m$ we have $|y_i|\leq m\beta_{i,i}$. Since
$Im(\tau_{g,g})=\beta_{g,g}- \beta^*_{g,g}$ and $\beta^*_{g,g}>0$
we also have $Im(\tau_{g,g})>\beta_{g,g}$. It follows from
\eqref{est1}  that for $n\in \Z^g$
\begin{multline*}
2Im(\tau)[n] -\,2\tra n y\,+2k Im(\tau_{g,g}) \geq
\sum_{i=1}^g \left(c \beta_{i,i}n_i^2-2|n_i|\beta_{i,i}m\right)+2k\beta_{g,g}\\
= \sum_{i=1}^g \left(c\beta_{i,i}\left(|n_i|-\frac{m}{c}\right)^2-\beta_{i,i}\frac{m^2}{c}\right)
+2k\beta_{g,g}\\
\geq  \sum_{i=1}^g \left(c d\left(|n_i|-\frac{m}{c}\right)^2-\beta_{g,g}\frac{m^2}{c}\right)
+2k\beta_{g,g}\\
=\sum_{i=1}^g c d\left(|n_i|-\frac{m}{c}\right)^2
+\left(2k\beta_{g,g}-g\beta_{g,g}\frac{m^2}{c}\right),
 \end{multline*}
and since $k>\frac{gm^2}{2c}$ we obtain
\[ 2Im(\tau)[n] -\,2\tra n y\,+k Im(\tau_{g,g}) > \sum_{i=1}^g c d\left(|n_i|-\frac{m}{c}\right)^2. \]
The above inequality together with \eqref{est2} imply
\begin{multline*}
|\theta(z,\tau)|\leq \sum_{n\in \Z} \eexp{-\pi \sum_{i=1}^g cd(|n_i|-m/c)^2}
\\=
\sum_{n\in \Z^g}\prod_{i=1}^g \eexp{-\pi cd(|n_i|-m/c)^2}\\
\leq \left(\sum_{n\in \Z} \eexp{-\pi cd(|n|-m/c)^2} \right)^g.
\end{multline*}
Since the series $\sum_{n\in \N} \eexp{-\pi cd(|n|-m/c)^2}$ converges like
$\int_{-\infty}^{+\infty}e^{-ax^2}\,dx$ and does not
depend on $(z,\tau)$, the
function $\theta(z,\tau)$ is bounded on $\CT(C_m\la \beta^*_0\ra^{\ell_{1,1}>d})$. \qed

\begin{prop}\vlabel{prop2} For all real numbers $m,d,R>0$
there is an open set $W\subseteq \C^g\times \C^n$
containing the set  $\CT_{\leq R}(C_m^{\ell_{1,1}\geq d})$
such that
the restriction of the function
$\vartheta(z,\tau)$ to $W$ is definable in $\Rae$.
\end{prop}
\proof
Recall that
\begin{multline*} \CT_{\leq R}(C_m^{\ell_{1,1}\geq d})=\\
\left\{ (z,\tau) \in \CT(C_m^{\ell_{1,1} \geq d}) \colon
\bigwedge_{i=1}^g \bigl| Re(z_i)\bigr|\leq R, \,
 \bigwedge_{i\leq i,j\leq g} \bigl| Re(\tau_{i,j})\bigr|\leq R \right\}.
\end{multline*}

Let $m,d,R>0$ be  real numbers.

Since the function $\vartheta(z,\tau)$ is $\Z^g$-periodic in $z$ and $(2\Z)^n$-periodic in $\tau$,
for any $k\in \N$ the function $(z,\tau)\mapsto \eexp{2\pi i k \tau_{g,g}}\vartheta(z,2\tau)$ is
$\Z^g\times\Z^n$-periodic.

Let $k\in \N$ and $\beta^*\in Int(\Mib_g)$  be as in Proposition \ref{prop1}. We are
going to apply Theorem \ref{defpseta} to the function
\[
\theta(z,\tau)= \eexp{2\pi i k \tau_{g,g}}\vartheta(z,2\tau). \]
The set $U=\C^g\times \CH_g$ is an open subset of $\C^g\times
\C^n$, and $\theta(z,\tau)$ is holomorphic and $\Z^g\times
\Z^n$-periodic on $U$. The integral cone we take is $C_m$ and the
linear function $\ell_{1,1}(y,\beta)=\beta_{1,1}$. For
$\beta_0^*=(0,\beta^*)$, by Proposition \ref{prop1}, $ \CT(C_m\la
\beta_0^*\ra^{\ell >d})\subseteq U$, and the function
$\theta(z,\tau)$ is bounded on $ \CT(C_m\la
\beta_0^*\ra^{\ell_{1,1} >d})$.  Also, since $\beta^*\in
Int(\Mib_g)$, it follows from the definition of $C_m$ that
$\beta_0^*\in Int(C_m)$.

Thus all assumptions of Theorem
\ref{defpseta} hold, and there is a definable  open $V\subseteq \C^g\times
\C^n$ containing   $\CT_{\leq R}(C_m\la \beta_0^*\ra)^{\ell_{1,1}\geq 2d}$ such that the restriction of
$\theta(z,\tau)$ to $V$ is definable in $\Rae$. Since $C_m$ is a cone and $\beta_0^*\in C_m$,
we have $C_m\la \beta_0^*\ra \subseteq C_m$, hence
$\CT_{\leq R}(C_m^{\ell_{1,1}\geq 2d}) \subseteq V$.

Intersecting $V$ with the set $\{ (z,\tau) \in \C^g\times \C^n \colon  |Re(\tau_{g,g})| <2R\}$
if needed, we can assume that $|Re(\tau_{g,g})|<2R$ on $V$.

It is easy to see that the restriction of the function $v\mapsto
\eexp{2\pi i k v}$ to the set $\{ v\in \C \colon |Re(v)| <2R\}$ is definable in the structure
$\Rae$, hence the restriction of the function $(z,\tau)\mapsto \vartheta(z,2\tau)$ to $V$
is also definable.

It is immediate that the restriction of the function
$\vartheta(z,\tau)$ to the set $W=\{ (z,\tau) \colon
(z,\frac{1}{2}\tau)\in V \}$ is definable in $\Rae$, and it is not
hard to see that $W$ contains the set $\CT_{\leq 1/2R}(C_{1/2m}^{\ell_{1,1}\geq 4d})$.

Since $m,d,R$ were arbitrary, the proposition follows. \qed

\begin{thm}\vlabel{def-theta}  For all $a,b\in \R^g$ and all real numbers $m,d,R>0$
there is an open set $U\subseteq \C^g\times \C^n$ containing
$\CT_{\leq R}\bigl(C_m^{\ell_{1,1}\geq d}\bigr)$
such that
the restriction of the function
$\teta a b(z,\tau)$ to $U$ is definable in $\Rae$.
\end{thm}
\proof We fix $a,b\in \R^g$, and $m,d,R>0$. Recall that
\[ \teta a b (z,\tau)=\eexp{\pi i(\tra a \tau  a + 2\, \tra a(z+b))}\vartheta(z+\tau a +b,\tau).
\]

For $u\in \R^g$ we denote by $\snorm u$ the sup-norm of $u$,
namely $\snorm u =\sup\{ |u_i| \colon i=1,\dotsc,g\}$.

\medskip

Let $f_1\colon \C^g\times \C^n\to \C $ be the map
\[ f_1(z,\tau)= \tra a \tau  a + 2\, \tra a(z+b).\]
If   $(z,\tau)\in \CT_{\leq R}\bigl(C_m^{\ell_{1,1}\geq d}\bigr)$ then
$\snorm {Re(z)}\leq R$ and  $|Re(\tau_{i,j})|\leq R$, hence
\begin{multline*}
|Re(f_1(z,\tau))|=\left|\, \sum_{1\leq i,j \leq g}Re(\tau_{i,j})a_ia_j  +2\sum_{i=1}^{g}a_i(Re(z_i)+b_i)
\right| \leq
\\
g^2 R\snorm{a}^2 +2g\snorm a(R+\snorm b).
\end{multline*}
Hence there is $R_1>0$ such that the set $V_1= \{ w\in \C \colon | Re(w) | <R_1 \}$
contains the image of $\CT_{\leq R}\bigl(C_m^{\ell_{1,1}\geq d}\bigr)$ under $f_1$.

It follows then  that the restriction of the function $w\mapsto
\eexp{\pi i w }$ to the set $V_1$ is definable in $\Rae$. Let
$U_1=f_1^{-1}(V_1)$. Clearly $U_1$ is  open in $\C^g\times \C^n$,
it contains $\CT_{\leq R}\bigl(C_m^{\ell_{1,1}\geq d}\bigr)$  and
the restriction of the function
\[ (z,\tau)\mapsto  \eexp{\pi i(\tra a \tau  a + 2\, \tra a(z+b))} \]
to $U_1$ is definable in $\Rae$.

\medskip

Let $f_2\colon \C^g\times \C^n\to \C^g$ be the map $f_2(z,\tau)=
z+\tau a + b$. If $(z,\tau)\in \CT_{\leq
R}\bigl(C_m^{\ell_{1,1}\geq d}\bigr)$ and $w=f_2(z,\tau)$ then
for $i=1,\dotsc g$ we have
\[ |Re(w_i)| \leq R + gR \|a\|_s +\|b\|_s,\]
and
\[ |Im(w_i)| \leq |Im(z_i)| + \sum_{j=1}^{g}|Im(\tau_{i,j})| \cdot \snorm a.\]
Applying Fact \ref{M-reduced}(b), and the definition of $C_m$ (see
\ref{note:notec}), we obtain
\[ |Im(w_i)| \leq (m+g\snorm a) Im(\tau_{i,i}). \]
Let $R_2= R + gR \|a\|_s +\|b\|_s$ and $m_2=m+g\snorm a$. It is immediate that the image
of $\CT_{\leq R}\bigl(C_m^{\ell_{1,1}\geq d}\bigr)$  under the map
$(z,\tau)\mapsto (f_2(z,\tau),\tau)$ is contained in the
set
$\CT_{\leq R_2}\bigl(C_{m_2}^{\ell_{1,1}\geq d}\bigr)$.
Applying Proposition \ref{prop2} we can find a definable open set
$U_2\subseteq \C^g\times C^n$ such that the restriction of the function
$\vartheta(z+\tau a +b,\tau)$ to $U_2$ is definable in $\Rae$. We can take $U=U_1\cap U_2$. \qed

\section{Definability of complex tori and their embeddings}

\begin{ntn} For a $g\times g$-matrix $M$  and $1\leq j \leq g$ we will denote by $M_{(j)}$ the
$j$-th column of $M$.
\end{ntn}

\begin{ntn}\vlabel{note:edtou} Let $D$ be a polarization type. For $\tau\in \CH_g$ we will denote
 by $E^D_\tau\subseteq \C^g$ the fundamental parallelogram of the lattice $\Lambda^D_\tau$:
\[ E^D_\tau=
\left\{ \sum_{i=1}^g t_i\tau_{(i)} + \sum_{j=1}^g s_jD_{(j)} \colon 0\leq t_i,s_j <
1 \right\}, \] and by $\bar E_{\tau}^D$ its closure in $\C^g$
\[ \bar E^D_\tau=
\left\{ \sum_{i=1}^g t_i\tau_{(i)} + \sum_{j=1}^g s_jD_{(j)} \colon 0\leq t_i,s_j \leq 1 \right\}. \]

\end{ntn}

Clearly $E^D_\tau$ contains a unique representative for every
$\Lambda^D_\tau$-coset, and we will identify the underlying set of $\pav D \tau$
with $E^D_\tau$. Using this identification it is not hard to equip each $\pav D
\tau$, uniformly in $\tau$, with a semi-algebraic $\C$-manifold structure in the
sense of \cite[Section 2.2]{mild} (see also Definition \ref{def-manifold} in
Appendix), and obtain a semi-algebraic family of complex-analytic manifolds $\{ \pav
D \tau \colon \tau \in \CH_q\}$.

\medskip

Our main goal in this section is to show that in the case $D_{1,1}\geq 2$, for a sufficiently
large set $F\subseteq \CH_g$  the restriction of $\varphi^D$ to the
set $\{ (z,\tau)\in \C^g\times \CH_g \colon \tau\in F,\, z\in \bar E^D_\tau \}$
is definable in $\Rae$, and in particlular the family of maps
\[ \varphi^D_\tau\colon \pav D \tau \to \PP^N(\C),\; \tau\in F\]
is definable in $\Rae$.

\medskip

As in the previous section we fix a positive integer $g$ and let $n=g(g+1)/2$.
We also fix a polarization type $D=Diag(d_1,\dotsc,d_g)$.

\begin{ntn}\vlabel{note:cxdv} Let $V$ be a subset of $\CH_g$.
\begin{enumerate}
\item  We will denote by $\CX^D(V)$ the following subset
of $\C^g\times \CH_g$:
\[ \CX^D(V)=\{ (z,\tau)\in \C^g\times \CH_g \colon \tau\in V,\, z\in \bar E^D_\tau \}. \]
\item For a real $K>0$ we will denote by  $\CX^D_{<K}(V)$ the following subset
of $\C^g\times \CH_g$:
\begin{multline*} \CX^D_{<K}(V)=\bigl\{ (z,\tau)\in \C^g\times \CH_g \colon \tau\in V,\,
 z=\tau r +D r'\\
\text{ for some $r,r'\in \R^g$ with } \snorm{r}< K, \snorm{r'}< K \bigr\}.
\end{multline*}

\end{enumerate}
\end{ntn}

\begin{rem}\vlabel{rem54}

Note that for $K>1$, we have $\CX^D(V)\subseteq \CX^D_{<K}(V)$.

It is not hard to see that for a closed subset $V\subseteq \CH_g$ the set
$\CX^D(V)$ is closed in $\C^g\times \CH_g$, and
if $V\subseteq \CH_g$
is open then for any $K>0$ the set $\CX^D_{<K}(V)$ is an open subset of $\C^g\times \CH_g$.

\end{rem}

\subsection{Definability of theta functions over the Siegel fundamental domain}
We first show that restrictions of theta functions to the set
$\CX^D(\fF_g)$ is definable in $\Rae$.

\begin{claim}\vlabel{theta0} Let $D=Diag(d_1,\dotsc,d_g)$ be a polarization type. For any real
  $K>0$ there are $d,m,R >0$ such that $\CX^D_{<K}(\fF_g)\subseteq
  \CT_{\leq R}\bigl(C_m^{\ell_{1,1}\geq d}\bigr)$.
\end{claim}
\proof Let $(z,\tau)\in \CX^D_{< K}(\fF_g)$. Then $\tau=\alpha+i\beta\in \fF_g$, and
$z=\tau r+D r'$ for some $r,r'\in \R^g$ with $\snorm r < K$, and $\snorm{r'}< K$. If
$x=Re(z)$ and $y=Im(z)$ hen $x=\alpha r+Dr'$ and $y=\beta r$, with $\beta\in \Mi_g$.

For each $i=1,\dotsc, g$, we have $|y_i|\leq
\sum_{j=1}^{g}|\beta_{i,j}|\,|r_j|$, hence, by
Fact~\ref{M-reduced}(b), $|y_i|< \frac{1}{2}g K\beta_{i,\,i}$.
Thus $(y,\beta)\in C_m$ for $m=\frac{1}{2}gK$, and, by Fact
\ref{Igusa15}(b), $(z,\tau)\in \CT(C_m^{\ell_{1,1}\geq d})$ for
$d=\frac{\sqrt{3}}{2}$.

\medskip

For $i=1,\dotsc,g$ we also have $|x_i|\leq \sum_{j=1}^g
|\alpha_{i,\,j}||r_j| +d_i|r_j|'$. Since $d_1\leq d_2 \leq \dotsc
\leq d_g$, and $\tau\in \fF_g$, it follows that $|x_i|<
\frac{1}{2}gK+d_gK$, and we can take $R=\frac{1}{2}gK+d_gK$. \qed

Applying Theorem \ref{def-theta} and using Remark \ref{rem54} we
obtain the following two corollaries.

\begin{cor}\vlabel{theta05}  For
  any  $a,b\in \R^g$ and $K>0$ there is an open set $U\subseteq \C^g\times \CH_g$ containing
$\CX^D_{<K}(\fF_g)$ such that
the
    restriction of the function $\teta a b (z,\tau)$ to $U$   is
    definable in $\Rae$.
\end{cor}

\begin{cor}\vlabel{theta1}
  For every $a,b\in \R^g$ there is an open set $U\subseteq \C^g\times
  \CH_g$ containing $\CX^D(\fF_g)$ such that  the restriction of $\teta a
  b(z,\tau) $ to $U$ is definable in $\Rae$.
\end{cor}

\subsection{Definability of the map $\varphi^D(z,\tau)$ }
Our next goal is to show that in the case when $d_1\geq 2$ and
$G\leqslant Sp(2g,\Z)$ has finite index then the restrictions of
the map $\varphi^D(z,\tau)$ to sets $\CX^D(\fF_g(G))$ is definable
in $\Rae$.

We need the following auxiliary Claim.

\begin{claim}\vlabel{theta2} Let $K>0$, and
$M=\begin{pmatrix} \alpha & \beta \\ \gamma &\delta \end{pmatrix} \in Sp(2g,\R)$.
  For every $a^1,b^1\in \R^g$ there is an open set $U\subseteq \C^g\times \CH_g$ containing the
set $\CX^D_{<K}(M\cdot \fF_g)$ such that
  the restriction of the function
  \[ f^{a_1}_{b_1}(z,\tau)= \eexp{ \,-\pi i (\tra z\gamma\,
    \tra{\bigl(\gamma (M^{-1}\cdot \tau)+\delta\bigr)}z) \,}\cdot \teta {a^1}{b^1}(z,\tau) \] to
$U$ is definable in $\Rae$.
\end{claim}
\proof
We will use the following fact that follows from the theta transformation formula
(see \cite[p. 221 and 8.6.1 on p. 227]{BL}
\begin{fact}[The theta transformation formula]
Assume $M,a^1,b^1$ are as in the above claim. Then there are
 $a,b\in \R^g$, and $k,k_1\in \C$ such that for all $z\in \C^g$ and $\tau_1\in \CH_g$ we have

\begin{multline*}  \teta{a^1}{b^1}(z ,M\cdot \tau_1)=\\
k \sqrt{det(\gamma\tau_1+\delta)}
\eexp{\pi i (k_1+\tra z\gamma\, \tra(\gamma\tau_1+\delta)z)}
\teta{a}{b}(\tra (\gamma\tau_1 +\delta)z,\tau_1).
\end{multline*}
\end{fact}

We apply the theta transformation formula with $\tau=M\cdot \tau_1$ and
obtain
\[
f^{a_1}_{b_1}(z,\tau)=
k \sqrt{det(\gamma \tau_1 +\delta)}
\eexp{\pi i k_1}
\teta{a}{b}(\tra (\gamma \tau_1 +\delta)z,\tau_1).
\]
Since $k$ and $k_1$ are constants, and the function $\tau\mapsto
\sqrt{det(\gamma\tau+\delta)}$ is semi-algebraic, it suffices  to
find an open $U\subseteq \C^g\times \CH_g$ containing  $\CX^D_{<K}(M\cdot\fF_g)$
such that  the restriction of the function
\[(z,\tau)\mapsto
\teta {a}{b}(\tra (\gamma(M^{-1}\cdot \tau) +\delta)z,M^{-1}\cdot \tau)\] to
 $U$  is definable in
$\Rae$.

\medskip

Clearly the map $g (z,\tau) = (\tra(\gamma(M^{-1}\cdot \tau)+\delta) z,
M^{-1}\cdot \tau)$ is semi-algebraic and holomorphic map. Thus  it is sufficient to show
that there is an open  set $V\subseteq \C^g\times \CH_g$
containing the image of $\CX^D_{<K}(M\cdot \fF_g)$ under $g$ such that  the restriction of the
function $\teta a b(z,\tau)$ to $V$ is definable in $\Rae$.

\medskip
Let $(z,\tau)\in \CX^D_{<K}(M\cdot \fF_g)$, and $(z_1,\tau_1)=g(z,\tau)$. It is immediate that
 $\tau_1\in \fF_g$, $\tau=M\cdot \tau_1$, and $z_1=\tra(\gamma\cdot \tau_1+\delta)z$.

 Since $(z,\tau)\in \CX^D_{<K}(M\cdot \fF_g)$, there  are $r,r'\in \R^g$, with $\snorm{r}< K$,
 $\snorm{r'} <K$, such that $z=\tau r+D r'$ .

Let $r_1,r_1'\in \R^g$ such that $z_1= \tau_1 r_1+D r_1'$. It follows from Fact
\ref{mainfact} by direct computations (see also Remark \ref{remain}) that
$$\left(\begin{array}{c}
  r_1\\
  r_1' \\
\end{array}\right)=\left(%
\begin{array}{cc}
  I & 0 \\
  0 & D^{-1} \\
\end{array}%
\right)\tra{M}\left(%
\begin{array}{cc}
  I & 0 \\
  0 & D \\
\end{array}%
\right)\left(\begin{array}{c}
  r\\
  r' \\
\end{array}\right)
$$ or equivalently
$\begin{pmatrix} r_1 \\ Dr_1'
\end{pmatrix}= \tra{M}\begin{pmatrix} r \\ Dr'
\end{pmatrix}$.

Thus $\snorm{r_1} < K\cdot \snorm{M}\cdot\snorm{D}$ and $\snorm{r_1'}< K\cdot
\snorm{M}\cdot\snorm{D}$. It follows that for $K'=K\cdot \snorm{M}\cdot \snorm{D}$,
the image of $\CX^D_{<K}(M\cdot \fF_g)$ under $g$ is contained in the set $\CX^D_{<
K'}(\fF_g)$. By Corollary \ref{theta05} there is an open set $V$ as needed. \qed

\begin{cor}\vlabel{Gdefcor} Fix a polarization type $D=diag(d_1,\ldots, d_g)$ and assume $d_1\geq 2$. Then
\begin{enumerate} \item For every
  $M\in Sp(2g,\R)$ there is an open $U\subseteq \C^g\times \CH_g$
  containing $\CX^D(M\cdot\fF_g)$ such that the restriction of $\varphi^D(z,\tau)$
to $U$ is definable in the structure $\Rae$.

\item Let $G<Sp(2g,\Z)$ be a subgroup of finite index and fix $\fF_g(G)$ a Siegel
fundamental set for the action of $G$ on $\CH_g$. Then there is an open set
$U\subseteq \C^g\times \CH_g$ containing $\CX^D(\fF_g(G))$ such that the restriction
of $\varphi^D(z,\tau)$ to $U$ is definable in the structure $\Rae$.
\end{enumerate}
\end{cor}
\proof  (1)
Recall that
\[ \varphi^D (z,\tau)= \left(\teta {c_0} 0(z,\tau)\,
:\, \teta {c_1} 0(z,\tau) :\,
 \dotsc\, :\,\teta {c_N} 0(z,\tau) \right). \]
Since each $f_a^b(z,\tau)$ is a product if $\teta a b (z,\tau)$ and a function
that does not depend on $a,b$, we have
\begin{multline*}
\left(\teta {c_0} 0(z,\tau)\,
:\, \teta {c_1} 0(z,\tau) :\,
 \dotsc\, :\,\teta {c_N} 0(z,\tau) \right)=\\
\left(f^{c_0}_0(z,\tau)\,
:\, f^{c_1}_0(z,\tau) :\,
 \dotsc\, :\,f^{c_N}_0(z,\tau) \right),
\end{multline*}
hence
$\varphi^D (z,\tau)=\left(f^{c_0}_0(z,\tau)\,
:\, f^{c_1}_0(z,\tau) :\,
 \dotsc\, :\,f^{c_N}_0(z,\tau) \right)$.
We now use Claim \ref{theta2} and Remark \ref{rem54} to get required set $U$.

To see (2), just note that $\fF_g(G)$ is given by finitely many $Sp(2g,\Z)$-translates of $\fF_g$
so we can apply (1).
\qed

\section{Projective embeddings of some classical families and moduli spaces}
In this section we use our  previous results to establish the definability
embeddings of some classical families of abelian varieties and their moduli spaces
into projective space.

\medskip

For the rest of this section we fix a polarization type $$D=Diag(d_1,\dotsc,d_g).$$

\subsection{Modular forms and the moduli space of principally
polarized  abelian varieties}

\begin{defn} Let  $\Gamma$ be a finite index subgroup of $Sp(2g,\Z)$, a
holomorphic $\phi:\CH_g\to \C$ is called {\em a modular form of
weight
$k$ with respect to $\Gamma$}, if for every $M=\left(%
\begin{array}{cc}
  \alpha & \beta \\
  \gamma & \delta \\
\end{array}%
\right)\in \Gamma$ and for every $\tau \in \CH_g$,
$$\phi(M\cdot \tau)=det(\gamma\tau+\delta)^k\phi(\tau).$$
\end{defn}

Clearly, if $\Gamma_1\leqslant \Gamma_2\leqslant Sp(2g,\Z)$ then every modular form
with respect to $\Gamma_2$ is also a modular form with respect to $\Gamma_1$. For
any $k\geq$ we let
$$\Gamma_g(k)= \{ M\in Sp(2g,\Z) \colon M \equiv I_{2g}\quad  \text{mod}(k)
\} $$

Using the definability of the theta functions we can now prove:
\begin{thm}\vlabel{modular-def} Let $\Gamma\leqslant Sp(2g,\Z)$ be a congruence
subgroup of $Sp(2g,\Z)$, namely a group which contains some $\Gamma_g(n)$. If $\phi$
is a modular form with respect to $\Gamma$ then there is a definable open $U\sub
\CH_g$ containing $\fF_g$ such that $\phi\rest U$ is definable in $\Rae$.
\end{thm}
\proof
 By Igusa's work \cite[Corollary, p. 235]{Igusa-paper} (see also \cite[Theorem III.6.2]{Bost}),
each such modular form with respect to the congruence subgroup $\Gamma_g(2k)$ is
algebraic over a ring generated by monomials in the so-called theta constants,
namely a ring generated by products
$\teta{a_1}{b_1}(0,\tau)\cdots\teta{a_k}{b_k}(0,\tau)$, for certain $a_i,b_i\in
\Q^g$. Since $\Gamma_g(2n)\leqslant \Gamma_g(n)\leqslant \Gamma$, it follows  that
every modular form with respect to $\Gamma$ is algebraic over the same ring. By
Corollary \ref{theta05}, the restriction of each $\teta{a}{b}(0,\tau)$ to some open
set containing $\fF_g$ is definable in $\Rae$. It is thus sufficient to note that
every function in the algebraic closure of the ring generated by the theta functions
is definable on some open set containing $\fF_g$.

 Assume that $h(\tau)$ is
the zero of a polynomial over the ring of theta constants. So
$h(\tau)$ is a zero of a polynomial
$$p_{\tau}(w)=\vartheta_n(\tau)w^n+\cdots+\vartheta_1(\tau)w+\vartheta_0(\tau),$$
where each $\vartheta_i$ is a polynomial in theta constants. We find a definable set
$U\sub \CH_g$ containing $\fF_g$ on which all $\vartheta_i$'s are definable and
consider the set
$$H=\{(\tau,w)\in U\times \C:p_{\tau}(w)=0\},$$ which is definable in $\Rae$.  By o-minimality we can write  $U$ as the union of finitely many
definable open simply connected sets and a set of smaller
dimension such that on each of these open sets the set $H$ is the
union of finitely many graphs of holomorphic functions. Now, on
each of the open sets the graph of $h$ equals one of the branches
of this cover, hence definable. The rest of the graph of $h$ is
obtained by taking topological closure. Thus, $h\rest U$ is
definable in $\Rae$.\qed

By the work of Baily and Borel, \cite[Theorem 10.11]{baily-borel}
(see also \cite[Theorem 7]{Geer}), the quotient
$Sp(2g,\Z)\backslash \CH_g$ can be embedded into projective space
using Siegel modular forms as coordinate functions. Namely, there
is a holomorphic map $F:\CH_g\to \PP^N(\C)$, given in coordinates
by modular forms with respect to $Sp(2g,\Z)$, which induces an
embedding of $Sp(2g,\Z)\backslash \CH_g$ into $\PP^N(\C)$.
Moreover, its image is Zariski open in some algebraic variety.
Taken together with Theorem \ref{modular-def} we can conclude the
definability in $\Rae$ of a restricted $F$.

\begin{thm} \vlabel{ppp} There is a definable open set $U\sub \CH_g$ containing $\fF_g$
such that $F\rest U$ is definable in $\Rae$.

More explicitly, the function $F$ satisfies:  for every
$\tau,\tau'\in U$, $F(\tau)=F(\tau')$ if and only if $\tau$ and
$\tau'$ are in the same $Sp(2g,\Z)$-orbit,  and  $F(U)$ is a
quasi-projective subset of $\PP^N(\C)$.
\end{thm}

\subsection{Embeddings of families} We consider here the uniform definability of the embedding of
abelian varieties into projective space.

The following is  immediate from Corollary \ref{Gdefcor}.
\begin{cor}\vlabel{Gdef}  Let $G<Sp(2g,\Z)$ be a subgroup of finite index and fix $\fF_g(G)\sub \CH_g$ a Siegel fundamental set
for $G$.
If $d_1\geq 2$ then the family of maps
\[ \{ \varphi^D_\tau\colon \pav D \tau \to \PP^N(\C)\colon \tau\in
\fF_g(G) \} \] is definable in the structure $\Rae$. In particular, when $d_1\geq 3$ (see Fact \ref{lef}),
 the family of projective abelian varieties $\{ \varphi^D(\pav
D \tau)\colon \tau\in \fF_g(G)\}$ is definable.
\end{cor}

As the next theorem shows, we can omit the restriction on $D$ and still obtain a uniformly definable
 family of emebddings of polarized abelian varieties of type $D$.

\begin{thm}\vlabel{Gdef2} For any polarization type $D$ there is in $\Rae$ a definable set $F\sub \CH_g$, containing a representative
for every $G_D$-orbit, and a definable embedding of the family $\{ \pav D \tau
\colon \tau\in F \}$ into projective space. More precisely, for some $N$ there is a
definable family of embeddings $\{h_{\tau}:\pav D \tau \to \PP^N(\C):\tau\in F\}$.
\end{thm}
\proof  If $d_1\geq 3$ then the theorem follows from Corollary
\ref{Gdef} with $h_{\tau}=\varphi^D_{\tau}$.

Assume $d_1<3$ and let $G=G_D\cap Sp(2g,\Z)$. It is known (see Fact \ref{0fin}
below) that $G$ has finite index in both groups, so we fix a Siegel fundamental set
$\fF_g(G)$ for $G$.
 For any $\tau\in \CH_g$ the map $z\mapsto 3z$ is an isomorphism from
  the torus $\CE^D_\tau$ onto the torus $\CE^{3D}_{3\tau}$, hence we have
in $\Rae$ a definable embedding of the family $\{\CE^D_\tau \colon \tau\in
\fF_g(G)\}$ into the family $\{\CE^{3D}_{3\tau} \colon \tau\in \fF_g(G)\}$.

Let $\gamma_1,\dotsc,\gamma_k\in Sp(2g,\Z)$ be such that $\fF_g(G)=\cup \gamma_i\cdot \fF_g$.
Let $\tau\in \gamma_i\cdot \fF_g$ and let $\tau_1= \gamma_i^{-1}\cdot \tau$.
Obviously $\tau_1\in \fF_g$ and $3\tau=\gamma_i\cdot(3\tau_1)$. Since $\tau_1\in \fF_g$,
then it is easy to see that $3\tau_1$ satisfies the conditions
  (i) and (iii) from the definition of $\fF_g$, and $3\tau_1 \in \fF_g+\beta$ for some
$\beta\in M_g(\Z)$ with $\snorm \beta \leq 2$. By the definition of the action of
$Sp(2g,\Z)$ on $\CH_g$, for any   $\beta\in M_g(\Z)$
we have $\fF_g+\beta=M_\beta\cdot \fF_g$, where $M_\beta=
\begin{pmatrix} I_g & \beta \\ 0 & I_g \end{pmatrix}\in Sp(2g,\Z)$.

Thus for every $i=1,\dotsc,k$ the family
$\{\CE^{3D}_{3\tau} \colon \tau\in \gamma_i\cdot \fF_g\}$ is contained in the family
$\{\CE^{3D}_\tau \colon \tau\in F_i\}$, where $F_i$ is the finite union
\[ F_i=\bigcup\bigl\{ \gamma_i\cdot M_\beta\cdot \fF_g \colon \beta\in M_g(\Z)
\text{ with } \snorm{\beta}\leq 2\bigr\}. \]
We can now use Corollary \ref{Gdefcor}(1) (and Fact \ref{lef}(b)) to uniformly embed each $\{\CE^{3D}_\tau \colon \tau\in F_i\}$
into a projective space, definably in the structure $\Rae$. \qed

\subsection{On embeddings  of  moduli spaces and universal families}

Note that as a corollary of Theorem \ref{Gdef2} we obtain, for every $g$ and $D$, a definable family
$\{h_{\tau}(\CE^D_{\tau}):\tau \in \fF_g(G_D)\}$ of
$g$-dimensional projective abelian varieties with polarization type $D$. Moreover, since $\tau$ varies
over a fundamental set for $G_D$, the family contains
 a representative for {\em every} $g$-dimensional polarized abelian variety of type $D$.
However, we do not claim that every such variety appears exactly once in the family. In order to
obtain such family we need to work with certain subgroups of $G_D$.

\begin{defn}[{see \cite[p. 233 (4) and p. 235, Lemma 8.9.1]{BL}}]\vlabel{GDD}
For a polarization type $D$
let $\GDD$ be the collection of all $2g\times 2g$ matrices
$\begin{pmatrix}
 \alpha & \beta\\
  \gamma& \delta \\
\end{pmatrix} $ in $G_D$ for which the following hold:
\begin{enumerate}[(i)]
\item There are $a,b,c,d\in M_g(\Z)$ such that $\alpha=I_g+Da$, $\beta=DbD$,
$\gamma=c$ and $\delta=I_g+dD$.
\item The diagonal elements of the matrices $D^{-1}\alpha\, \tra{\beta}D^{-1}$ and
$\gamma\, \tra{\delta}$ are even integers.
\end{enumerate}
\end{defn}
By (i), $\GDD< Sp(2g,\Q)\cap M_g(\Z)=Sp(2g,\Z)$.

\medskip

We recall the following.
\begin{fact}\vlabel{0fin} $\GDD$ has finite index in both $G_D$ and $Sp(2g,\Z)$.
\end{fact}
\proof By \cite[Lemma 8.9.1 (b)]{BL} $\GDD$ has finite index in $G_D$. It is not
hard to check by direct computations that $\GDD$ contains the principal congruence
subgroup $\Gamma_g(2(d_g)^2)$.  Since $\Gamma_g(2(d_g)^2)$ has finite index in
$Sp(2g,\Z)$, $\GDD$ has finite index as well. \qed

\medskip

Let $\Psi^D(\tau)\colon \CH_g \to \PP^N(\C)$ be the map
\[\Psi^D(\tau)= \varphi^D(0,\tau),\] and
 $\Phi^D\colon \C^g\times \CH_g\to \PP^N(\C)\times \PP^N(\C)$ be the map
\[ \Phi^D(z,\tau)=(\varphi^D(z,\tau),\Psi^D(\tau)). \]
\vlabel{note:phipsi}

We have a commuting diagram
\[
\begin{diagram}
\node{\C^g\times \CH_g}\arrow{e,t}{\Phi^D}\arrow{s,l}{\pi_2} \node{\PP^N(\C)\times \PP^N(\C)}
 \arrow{s,r}{\pi_2}\\
\node{\CH_g}\arrow{e,t}{\Psi^D} \node{\PP^N(\C)}
\end{diagram}
\]

\medskip

The reason for introducing $G_D(D)_0$ is the following result, which can be seen as
also saying that for appropriate $D$, the map $\Psi^D$ induces an embedding of the
quotient space $G_D(D)_0 \backslash \CH_g$ into projective space, whose image is a
quasi-projective variety. For a proof we refer to \cite[Theorem 8.10.1. and Remark
8.10.4]{BL}.
\begin{fact}\vlabel{fact-moduli}  Assume $d_1\geq 4$ and $2|d_1$ or $3|d_1$.
The map $\Psi^D\colon \CH_g \to\PP^N$ is an immersion,
the image of $\CH_g$ under $\Psi^D$ is a Zariski open subset of an algebraic variety in
$\PP^{N}$, and $\Psi^D(\tau)=\Psi^D(\tau')$ if and only if
there exists $M\in G_D(D)_0$ with $M\cdot \tau=\tau'$.

Moreover,
if $\Psi^D(\tau)=\Psi^D(\tau')$ then $\varphi^D_{\tau}(\pav D \tau)=\varphi^D_{\tau'}(\pav D {\tau'})$.
\end{fact}

\begin{ntn}\vlabel{note:fFgd} We will denote by $\fF_g^D$ a fixed Siegel  fundamental set  for the action $\GDD$ on $\CH_g$. \end{ntn}
Since $\GDD$ is a subgroup of $G_D$, by Fact \ref{mainfact}, every polarized abelian variety of type $D$
 is isomorphic to one of $\pav D \tau$ with $\tau\in \fF_g^D$.
Note also that by Fact \ref{fact-moduli}, for $D$  as above, if
$\fF\sub \CH_g$ contains $\fF_g^D$ then
$\Psi^D(\fF)=\Psi^D(\CH_g)$ and
$\Phi^D(\CX^D(\fF))=\Phi^D(\C^g\times \CH_g)$. An immediate
corollary of Corollary \ref{Gdefcor}(2) is:

\begin{thm}\vlabel{moduli2} If $d_1\geq 4$ and $2|d_1$ or $3|d_1$ then there is a definable open set $U\sub \C^g\times \CH_g$
containing $\CX^D(\fF^D_g)$ such that $\Phi^D\rest U$ is definable
in $\Rae$.
\end{thm}

\subsection{A new proof for a theorem of Baily}
 In this section we demonstrate
how o-minimality can be used to provide an alternative proof of
the following theorem of Baily, \cite{baily}.

\begin{thm}[Baily]\vlabel{baily}  Assume $d_1\geq 4$ and
$2|d_1$ or $3|d_1$. Then the  image of $\C^g\times \CH_g$ under $\Phi^D$ is a Zariski open subset
of an algebraic variety in $\PP^N(\C)\times \PP^N(\C)$.
\end{thm}
\proof Let $U\sub \C^g\times \CH_g$ be the open set from Theorem
\ref{moduli2}.  The map $\Phi\rest U$ is definable in $\Rae$ and
$\Phi(U)=\Phi(\C^g\times \CH_g)$.

Let $X=\Phi^D(\CX^D(\CH_g^D))$. The main step in proving the theorem is to show that
the topological closure $\bar X$ of $X$ in $\PP^N(\C)$ is a projective variety.
For that we will use the following theorem, proved in Appendix.
For $Z\subset \C^k$, we let $Fr(Z)=Fr_{\C^k}(Z)=Cl(Z)\setminus Z$,
where $Cl(Z)$ stands for the topological closure of $Z$ in $\C^k$.
We write $\dim_{\R}(Z)$ for the o-minimal dimension of $Z$.

\begin{fact}\vlabel{compat}
 Let $\phi:U\to \PP^N(\C)\times \PP^N(\C)$
be a  finite-to-one holomorphic map from an open $U\sub \C^m$ which is definable in
$\Rae$. Assume that there is a definable set $F\sub U$   such that $\phi(U)=\phi(F)$
and such that $\dim_{\R}(Fr(F))\leq 2m-2$.
 Then the topological closure of $\phi(U)$ in $\PP^N(\C)\times \PP^N(\C)$  is
an algebraic variety.
\end{fact}

We would like to  apply the above result to the map $\phi=\Phi^D\rest U$ with  $F=\CX^D (\fF^D_g)$ and $U$ viewed
as an open subset of $\C^g\times \C^n$.
Since the pre-image of every point under the map
$\Phi^D$ is a discrete subset  of $\C^g\times \CH_g$, by o-minimality its intersection with $U$
must be finite, so $\phi^D$ is finite-to-one on $U$. It is left to see that $\dim_{\R}(Fr_{\C^g\times \C^n}(F))\leq 2(g+n)-2$.

For  every $\gamma\in Sp(2g,\R)$, we consider the definable map $\tau\mapsto \gamma\tau$, from the open set
$\CH_g\sub \C^n$ into $\C^n$. By Fact \ref{Igusa15}, the set $\fF_g$
a closed subset of $\C^n$ and therefore, by Theorem \ref{complexman} below,
$\dim_{\R}(Fr_{\C^n}(\gamma \cdot \fF_g))\leq 2n-2$.  Since $\fF^D_g$ is a finite union of  such translates (see (\ref{fFD})), it follows that
 $\dim_{\R}(Fr_{\C^n}(\fF^D_g))\leq 2n-2$. Recall that
 $F=\CX^D (\fF^D_g)$ is defined as $$\left\{(\sum_{i=1}^gt_i\tau_{(i)}\! +\! \sum_{j=1}^gs_jD_{(j)},\tau)\in \C^g\times \CH_g: \tau \in \fF^D_g
 \, , \, 0\leq t_i,s_j\leq 1\right\}.$$

 Consider the frontier of $F$ inside $\C^g\times \C^n$. It is easy to see that if $(z,\tau)$ belongs
 to this frontier then $\tau\in Fr_{\C^n}(\fF^D_g)$, so by the above,  the
 real dimension
 of $Fr(F)$ is at most  $2g+2n-2$, as we wanted.

 We can therefore apply Fact \ref{compat} and conclude that $\bar X$ is an algebraic subvariety of $\PP^N(\C)\times \PP^N(\C)$.

Let $Y=\pi_2(X)$ and $\bar Y=Y$ be the topological closure of $Y$. By Fact \ref{fact-moduli},
$\bar Y$ is a Zariski closed set and $Y$ is Zariski open in $\bar Y$.

It is not hard to see that $Fr(X)$ is disjoint from
$\pi^{-1}_2(Y)$, hence $X=\bar X\cap \pi^{-1}_2(Y)$. Since $Y$ is
Zariski open in $\bar Y$, the set $X$ is Zariski open in $\bar X$.
\qed

\section{Appendix: O-minimality and complex analysis}

We review here some basic notions from \cite{mild} and \cite{ICM} and prove the results we used earlier.

\begin{defn} \label{def-manifold} Let $\CR=\la \R;,+,\cdot, \cdots\ra$ be an o-minimal expansion of the real field.
{\em A definable $n$-dimensional $\C$-manifold} is a definable set $M$, equipped with a finite cover of definable sets
$M=\bigcup_i U_i$, each in definable bijection with an open subset of $\C^n$, and  such that the transition maps are
holomorphic.
\end{defn}

For $M$ a manifold and $X\sub M$ we write $Cl_M(X)$ for the topological closure of $X$ in $M$, and
 $Fr_M(X)$ for the frontier $Cl_M(X)\setminus X$.
We fix $\CR$ an o-minimal expansion of the field of reals. All definability is in $\CR$.

\begin{thm} \vlabel{complexman}
Assume that $\phi:U\to N$ is a definable
holomorphic map from a definable open subset of $\C^n$ into a definable $\C$-manifold.
\begin{enumerate}
\item  Let $F\sub U$ be a closed (not necessarily definable) set. Then
there is a definable $Y\sub N$, with $\dim_{\R}(Y)\leq 2n-2$ such that $Fr_N(\phi(F))\sub Y$. Moreover, the same result holds
if we assume, instead of $F$ being closed, that $Fr_{\C^n}(F)$ is contained in some definable set
whose real dimension is $2n-2$.

\item Assume that $\phi$ is a finite-to-one map, and that there is a set $F\sub U$ with $\phi(F)=\phi(U)$,
such that $ Fr_{\C^n}(F)$ is contained in a definable set of dimension $2n-2$.
 Then $Cl_N(\phi(U))$ is a complex analytic
subset of $N$. In particular, if $N=\PP^k(\C)$ then, by Chow's Theorem, $Cl_N\phi(U))$ is a projective variety.
\end{enumerate}
\end{thm}
\proof  (1) Assume first that $F$ is closed. Note that if $F$ is compact then
$\phi(F)$ is closed so $Fr_N(\phi(F))$ is empty. Hence, the points of
$Fr_N(\phi(F))$ arise from the behavior of $\phi$ ``at $\infty$''.

We repeat the argument in \cite[Section 7.2]{ICM}. We write $\PP^n(\C)=\C^n\cup H$,
for $H$ a projective hyperplane, and let $\Gamma\sub \PP^n(\C)\times N$ be the
closure of the graph of $\phi$. We let $\pi:\PP^n(\C)\times N\to N$ be the
projection onto the second coordinate. Because $F$ is closed, the frontier of
$\phi(F)$ in $N$ is contained $Y=\pi(\Gamma\cap (H\times N))$, namely for every
$y\in Fr_N(\phi(F))$ there exists $z\in H$ such that $(z,y)\in \Gamma$. It is
therefore sufficient to see that $\dim_{\R}(\Gamma\cap (H\times N))\leq 2n-2$.

We let $B_{inf}$ be the set of all $(z,y)\in \Gamma\cap (H\times N)$ such that there
are infinitely $y'\in N$ with $(z,y')\in \Gamma$ and let $B_{fin}=\Gamma\cap
(H\times N)\setminus B_{inf}$. By \cite[Lemma 6.7(ii)]{crelle},
$\dim_{\R}(B_{inf})\leq 2n-2$. But now, since $\dim_{\R}(H)=2n-2$ it follows from
the definition of $B_{fin}$ that also $\dim_{\R}(B_{fin})\leq 2n-2$. Since
$\Gamma\cap (H\times N)\sub B_{inf}\cup B_{fin}$, we have $\dim_{\R}(\Gamma\cap
(H\times N))\leq 2n-2$, as required.

As for the ``moreover'' statement, we can repeat the above argument, with $H$ replaced by $H'=Fr_{\PP^n(\C)}(F)$.
The assumption implies $H'$ is contained in a definable subset of $\PP^n(\C)$ of dimension $2n-2$, so we can repeat the
argument.

(2)
By \cite[Corollary 6.3]{crelle}, there is a definable closed $E\sub N$ with
$\dim_{\R}(E)\leq \dim_{\R} \phi(U)-2$, such that $A=\phi(U)\setminus E$ is locally $\C$-analytic in $N$.
Since $\phi$ is finite-to-one, the intersection of  $\phi(U)$ with every open subset of $N$ is either empty
or of dimension $2n-2$. Because $\phi(U)=\phi(F)$, we can conclude from (1) then
 $\dim_{\R}Fr_N(\phi(U))\leq 2n-2$. In particular, $\dim_{\CR}Fr_N(A)\leq 2n-2$.

 We can now apply \cite[Theorem 4.1]{crelle} to $A$ and conclude that $Cl(A)$ is an analytic subset of $N$.
  It is easy to see that $Cl_N(A)=Cl_N(\phi(U))$.\qed

\end{document}